%% file: Phase_retrieval_arxiv.tex
\documentclass[a4paper,11pt]{article}
\usepackage{authblk}
\usepackage[british]{babel}
\usepackage[mono=false]{libertine}
\usepackage{setspace}

\usepackage[utf8]{inputenc} 
\usepackage[T1]{fontenc}    
\usepackage{url}            
\usepackage{booktabs}       
\usepackage{amsfonts}       
\usepackage{nicefrac}       
\usepackage{microtype}      
\usepackage{mathtools}
\usepackage{amsmath,amsfonts,amssymb,amsthm,mathrsfs,bbm}
\usepackage{latexsym,amscd,amsbsy,dsfont}
\usepackage{mathtools}
\usepackage{dsfont}
\usepackage{color}
\usepackage[usenames, dvipsnames]{xcolor}
\usepackage{graphicx}
\usepackage{epstopdf}
\usepackage{float}
\usepackage{cases}
\usepackage[top=2.5cm, bottom=2.5cm, left=2cm, right=2cm]{geometry}
\usepackage{caption}
\usepackage{physics}
\usepackage{tensor}
\usepackage{subcaption}
\usepackage[most]{tcolorbox}
\tcbuselibrary{theorems}
\usepackage{empheq}

\usepackage{xr-hyper} 
\usepackage[pdfpagemode=UseNone,bookmarksopen=false,colorlinks=true,urlcolor=RoyalBlue,citecolor=YellowOrange,citebordercolor=RoyalBlue,linkcolor=RoyalBlue]{hyperref}
\usepackage{bookmark}

\usepackage{times}
\usepackage{enumitem}

\author{Benjamin Aubin$^\dagger$, Bruno Loureiro$^\dagger$, Antoine Baker$^\star$ \\ 
Florent Krzakala$^\star$, Lenka Zdeborov{\'a}$^\dagger$}
\date{
$\dagger$ \textit{Institut de Physique Th\'eorique \\
CNRS \& CEA \& Universit\'e Paris-Saclay, Saclay, France}\\
$\star$ \textit{Laboratoire de Physique Statistique\\
CNRS \& Sorbonnes Universit\'es \& \\
\'Ecole Normale Sup\'erieure, PSL University, Paris, France}\\
\vspace{1cm}
}

\title{Exact asymptotics for phase retrieval and compressed sensing with random generative priors}
\input{settings_custom.tex}

\begin{document}
\maketitle

\begin{abstract}
We consider the problem of compressed sensing and of (real-valued) phase retrieval with random measurement matrix.
We derive sharp asymptotics for the information-theoretically optimal performance and for the best known polynomial algorithm for an ensemble of generative priors consisting of fully connected deep neural networks with random weight matrices and arbitrary activations. 
We compare the performance to sparse separable priors and conclude that generative priors might be advantageous in terms of algorithmic performance. 
In particular, while sparsity does not allow to perform compressive phase retrieval efficiently close to its information-theoretic limit, it is found that under the random generative prior compressed phase retrieval becomes tractable. 
\end{abstract}


\input{files/introduction.tex}
\input{files/bayesian.tex}

\input{files/phase_diagrams.tex}

\input{files/conclusion.tex}

\section*{Acknowledgements}
This work is supported by the ERC under the European Union’s Horizon
2020 Research and Innovation Program 714608-SMiLe, as well as by the
French Agence Nationale de la Recherche under grant
ANR-17-CE23-0023-01 PAIL. 
We gratefully acknowledge the support of
NVIDIA Corporation with the donation of the Titan Xp GPU used for this
research. We thank Google Cloud for providing us access to their platform through the Research Credits Application program. BL was partially financed by the Coordena\c{c}\~ao de Aperfei\c{c}oamento de Pessoal de N\'ivel Superior - Brasil (CAPES) - Finance Code 001.

We thank Antoine Maillard for insightful discussions in the early stages of this work and Andr\'e Manoel for sharing his MLAMP code with us.

\clearpage
\bibliographystyle{unsrt}
\bibliography{refs}


\end{document}

%% file: settings_custom.tex
\def \({\left(}
\def \){\right)}
\def \[{\left[}
\def \]{\right]}

\newcommand{\bx}{{\mathbf {x}}}

\newcommand{\be}{\begin{equation}}
\newcommand{\ee}{\end{equation}}
\newcommand{\beqa}{\begin{eqnarray}}
\newcommand{\eeqa}{\end{eqnarray}}

\newcommand{\bea}{\begin{align}}
\newcommand{\eea}{\end{align}}

\DeclareMathAlphabet{\varmathbb}{U}{bbold}{m}{n}

\newcommand{\bbR}{\mathbb{R}}

\newcommand{\mO}{\mathcal{O}}

\newcommand{\extr}{\textrm{\textbf{extr}}}
\newcommand{\td}[1]{{\tilde{#1}}}

\renewcommand{\subparagraph}[1]{\vspace{0.5cm} $\bullet$ \textit{#1} \\ \vspace{0.2cm} }

\renewcommand{\vec}[1]{{\mathbf{#1}}}

\def\out{\text{out}}

\def\dd{\text{d}}

\def\iid{\text{i.i.d.}}
\def\IT{\text{IT}}
\def\alg{\text{alg}}

\def\extr{{  \textbf{extr}}}

\renewcommand{\vec}[1]{{\textbf{#1}}}
\newcommand{\mat}[1]{{\textrm{#1}}}

%% file: files/introduction.tex
\section{Introduction}
\label{sec:intro}
Over the past decade the study of compressed sensing has lead to  significant developments in the field of signal processing, with novel sub-Nyquist sampling
strategies and a veritable explosion of work in sparse representation. A central observation is that sparsity allows one to measure the signal with fewer observations than its dimension \cite{donoho2006compressed,candes2006near}.
The success of neural networks in the recent years suggests another powerful and generic way of representing signals with multi-layer generative priors, such as those used in generative
adversarial networks (GANs) \cite{goodfellow2014generative}. It is therefore natural to replace sparsity by generative neural network models in compressed sensing and other inverse problems, a strategy that was successfully explored in a number of papers, e.g. \cite{tramel2016approximate,tramel2016inferring,bora2017compressed,manoel2017multi,hand2017global,fletcher2018inference,hand2018phase,mixon2018sunlayer,aubin2019spiked}.
While this direction of research seems to have many promising applications, a systematic theory of what can be efficiently achieved still falls short of the one developed over the past decade for sparse signal processing. 
Our aim is therefore to dialogue with the broad program of studying how generative models can help solving inverse problems using the toolbox of statistical physics.
In this manuscript, we build on a line of work allowing for theoretical analysis in the case the measurement and the weight matrices of the prior are random \cite{manoel2017multi,reeves2017additivity,fletcher2018inference,gabrie2018entropy,aubin2019spiked}. 

We employ tools originally developed in the context of statistical physics to derive precise asymptotics for the information-theoretically optimal thresholds for signal recovery and for the performance of the best known polynomial algorithm in two such inverse problems: (real-valued) phase retrieval and compressed sensing. These two problems of interest can be framed as a \emph{generalised linear estimation}. Given a set of observations $\vec{y}\in\mathbb{R}^{n}$ generated
from a fixed (but unknown) signal $\vec{x}^{\star}\in\mathbb{R}^{d}$ as 
\begin{equation}
\vec{y}=\varphi\left(\rm{A}\vec{x}^{\star}\right),
\end{equation}
the goal is to reconstruct $\vec{x}^{\star}$ from the knowledge of $\vec{y}$, $\varphi$ and $\mat{A}\in\bbR^{n \times d}$. Compressed sensing and phase retrieval are particular instances of this problem, corresponding to $\varphi(x)= x$ and $\varphi(x) = |x|$ respectively. Two key questions in these inverse problems are a) how many observations $n$ are required for theoretically reconstructing the signal $\vec{x}^{\star}$, and b) how this can be done in practice - i.e. to find an efficient algorithm for reconstruction. Signal structure plays an important role in the answer to both these questions, and have been the subject of intense investigation in the literature. A typical situation is to consider signals admitting a low-dimensional representation, such as sparse signals, for which $k-d$ of the $d$ components of $\vec{x}^{*}$ are exactly zero, see e.g. \cite{candes2015phase, netrapalli2013phase}.

In this manuscript, we consider instead structured signals drawn from a generative model $\vec{x}^{\star} = G(\vec{z})$, where $\vec{z}\in\mathbb{R}^{k}$ is a low-dimensional latent representation of $\vec{x}^{\star}$. In particular, we will focus in generative multi-layer neural networks, and in order to provide a sharp asymptotic theory, we will restrict the analysis to an ensamble of random networks with known random weights:
\begin{equation}
\vec{x}^{\star} = G\left(\vec{z}\right) = \sigma^{(L)}\left(\mat{W}^{(L)}\sigma^{(L-1)}\left(\mat{W}^{(L-1)}\cdots\sigma^{(1)}\left(\mat{W}^{(1)}\vec{z}\right)\cdots\right)\right), 
\label{eq:generative_model}
\end{equation}
\noindent where $\sigma^{(l)}:\mathbb{R}\to\mathbb{R}$, $1\leq l\leq L$ are component-wise non-linearities. As aforementioned, we take ${\rm{A}}\in\mathbb{R}^{n\times d}$ and ${\mat{W}}^{(l)}\in\mathbb{R}^{k_{l}\times k_{l-1}}$ to have i.i.d. Gaussian entries with zero means and variances $1/d$ and $1/k_{l-1}$ 
respectively, and focus on the high-dimensional regime defined by taking $n, d, k_{l} \to \infty$ while keeping the measurement rate $\alpha = n/d$ and the layer-wise aspect ratios $\beta_{l} = k_{l+1}/k_{l}$ constant. We stress that in this regime the depth $L$ is of order one when compared to the width of the generative network, which scales with the input dimension $d$. With this observation in mind, we adopt the standard terminology in machine learning of denoting networks with $L>1$ as \emph{deep}.
To provide a comparison with previous results for sparse signals, it is useful to define the total compression factor $\rho = k/d$. We note, however, that the comparison between generative and sparse priors herein is not based on a quantitative comparison between the reconstruction estimation errors. Indeed, since data is generated differently in both cases, such a comparison would make little sense. Instead, we compare qualitative properties of the phase diagrams, taking as a surrogate for algorithmic hardness the size of the statistical-to-algorithmic gap in these two different reconstruction problems.
Our results hold for latent variables drawn from an arbitrary separable distribution $\vec{z}\sim P_{z}$, and for arbitrary activations $\sigma^{(l)}$, although for concreteness we present results for $\vec{z}\sim \mathcal{N}(\vec{0},\mat{I}_{k})$ and $\sigma^{(l)}\in\{\text{linear}, \text{ReLU}\}$, as it is commonly the case in practice with GANs. 

\paragraph{Previous results on sparsity:} Sparsity is probably the most widely studied type of signal structure in linear estimation and phase retrieval. It is thus instructive to recall the main results for sparse signal reconstruction in these inverse problems in the high-dimensional regime with random measurement matrices studied in this manuscript. Optimal statistical and algorithmic thresholds have been established non-rigorously using the replica-method in a series of works \cite{wu2012optimal,krzakala_statistical-physics-based_2012,reeves2012compressed,zdeborova2016statistical}. Later the information theoretic results, as well as the corresponding minimum mean squared error (MMSE), has been rigorously proven in \cite{barbier2016mutual,Reeves,Barbier2017c}. So far, the best known polynomial time algorithm in this context is the \emph{Approximate Message Passing} (AMP) algorithm, the new avatar of the mean-field approach pioneered in statistical mechanics \cite{mezard1987spin}, that has been introduced in \cite{donoho2009message,rangan2011generalized, krzakala_statistical-physics-based_2012,schniter2014compressive,metzler2017coherent} for these problems, and can be rigorously analysed \cite{bayati2011dynamics}. For both (noiseless) compressed sensing and phase retrieval, the information theoretic limit for a perfect signal recovery is given by $\alpha>\alpha_{\rm IT}=\rho_s$, with $\rho_s$ being the fraction of non-zero components of the signal~$\vec{x}^{\star}$.

The ability of AMP to exactly reconstruct the signal, however, is different. A non-trivial line $\alpha^{\rm sparse}_{\rm alg}(\rho_s)>\alpha_{\rm IT}$ appears below which AMP fails. No polynomial algorithm achieving better performance for these problems is known. Strikingly, as discussed in \cite{Barbier2017c}, the behaviour of the sparse linear estimation and phase retrieval is drastically different: while $\alpha^{\rm sparse}_{\rm alg}(\rho_s)$ is going to zero as $\rho_s\to 0$ for sparse linear estimation hence allowing for compressed sensing, it is not the case for the phase retrieval, for which  $\alpha^{\rm sparse}_{\rm alg} \to 1/2$ as  $\rho_s\to 0$. As a consequence, \emph{no efficient approach to real-valued compressed phase retrieval with small but order one $\rho_s$ in the high-dimensional limit is known}. 

\paragraph{Summary of results:}
In this work, we replace the sparse prior by the multi-layer generative model introduced in eq.~\eqref{eq:generative_model}. Our main contribution is specifying the interplay between the number of measurements needed for exact reconstruction of the signal, parametrised by $\alpha$, and its latent dimension~$k$. Of particular interest is the comparison between a sparse and separable signal (having a fraction $\rho_s$ of non-zero components) and the structured generative model above, parametrised by $\rho=k/d$. While the number of unknown latent variables is the same in both cases if $\rho=\rho_s$, the upshot is that generative models offer algorithmic advantages over sparsity. More precisely:
\begin{enumerate}[label=(\alph*)]
	\item We analyse the minimum mean square error (MMSE) of the optimal Bayesian estimator for the compressed sensing and phase retrieval problems with generative priors of arbitrary depth, choice of activation and prior distribution for the latent variable. We derive sufficient conditions for the existence of an \emph{undetectable phase} in which better-than-random estimation of $\vec{x}^{\star}$ is impossible, and characterise in full generality the threshold $\alpha_{c}$ beyond which partial signal recovery becomes statistically possible. 
	\item Fixing our attention on the natural choices of activations $\sigma\in\{\text{linear},\text{ReLU}\}$, we establish the threshold $\alpha_{\IT}$ above which perfect signal reconstruction is theoretically possible. This threshold can be intuitively understood with a simple counting argument.
	\item We analyse the performance of the associated \emph{Approximate Message Passing} algorithm \cite{manoel2017multi}, conjectured to be the best known polynomial time algorithm in this setting. This allows us to establish the algorithmic threshold $\alpha_{\alg}$ below which no known algorithm is able to perfectly reconstruct $\vec{x}^{\star}$. 	 
\end{enumerate}
As expected, the thresholds $\{\alpha_{c}, \alpha_{\IT}, \alpha_{\alg}\}$ are functions of the compression factor $\rho$, the number of layers $L$, the aspect ratios $\{\beta_{l}\}_{l=1}^L$ and the activation functions. In particular, for a fixed architecture we find that the algorithmic gap $\Delta_{\alg} = \alpha_{\alg} - \alpha_{\IT}$ is drastically reduced with the depth $L$ of the generative model, beating the algorithmic hindrance identified in \cite{Barbier2017c} for compressive phase retrieval with sparse encoding. 

%% file: files/bayesian.tex
\section{Information theoretical analysis}
\label{sec:it_analysis}

\subsection{Performance of the Bayes-optimal estimator}
In our analysis we assume that the model generating the observations $\vec{y}\in\mathbb{R}^{n}$ is known. Therefore, the optimal estimator minimising the mean-squared-error in our setting is given by the Bayesian estimator
\begin{align}
\hat{\vec{x}}^{\text{opt}} = \underset{\hat{\vec{x}}}{\text{argmin}}||\hat{\vec{x}}-\vec{x}^{\star}||^2_2 = \mathbb{E}_{P(\vec{x}|\vec{y})}\[ \vec{x}\]\,.
\label{eq:optimal_estimator}
\end{align}
The posterior distribution of the signal given the observations is in general given by:
\begin{align}
	P(\vec{x}|\vec{y}) &= \frac{1}{\mathcal{Z}(\vec{y})}P_{x}(\vec{x})\prod\limits_{\mu=1}^{n}\delta\left(y^{\mu}-\varphi\left(\sum\limits_{j=1}^{d} {\mat{A}}^{\mu}_{j}x_j\right)\right),
\end{align}
\noindent where the normalisation $\mathcal{Z}(\vec{y})$ is known as the \emph{partition function}, and $\varphi$ is the nonlinearity defining the estimation problem, e.g. $\varphi(x) = |x|$ for phase retrieval and $\varphi(x) = x$ for linear estimation. We note that the presented approach generalizes straightforwardly to account for the presence of noise, but we focus in this paper on the analysis of the noiseless case. For the generative model in eq.~\eqref{eq:generative_model}, the prior distribution $P_{x}$ reads
\begin{align}
P_{x}(\vec{x}) = \int_{\mathbb{R}^{k}}\dd \vec{z}~ P_{z}(\vec{z})\prod\limits_{l=1}^{L}\int_{\mathbb{R}^{k_{l}}}\dd \vec{h}^{(l)}~P^{(l)}_{\out}\left(\vec{h}^{(l+1)}\Big|\mat{W}^{(l)}\vec{h}^{(l)}\right)\,,
\label{eq:Px}
\end{align}
\noindent where for notational convenience we denoted $\vec{x} \equiv \vec{h}^{(L+1)}$, $\vec{z} \equiv \vec{h}^{(1)}$ and defined the likelihoods $P^{(l)}_{\out}$ parametrising the output distribution of each layer given its input. As before, this Bayesian treatment also accounts for stochastic activation functions, even though we focus here on deterministic ones.

Although exact sampling from the posterior is intractable in the high-dimensional regime, it is still possible to track the behaviour of the minimum-mean-squared-error estimator as a function of the model parameters. Our main results are based on the line of works comparing, on one hand, the information-theoretically best possible reconstruction, analysing the ideal Bayesian inference decoder, regardless of the computation cost, and on the other, the best reconstruction using the most efficient known polynomial algorithm - the approximate message passing.

Our analysis builds upon the statistical physics inspired multi-layer formalism introduced in \cite{manoel2017multi}, who showed using the cavity and replica methods that the minimum mean-squared-error achieved by the Bayes-optimal estimator defined in eq.~\eqref{eq:optimal_estimator} can be written, in the limit of $n,d \to \infty$ and $\alpha=n/d=\Theta(1)$ for a generic prior distribution $P_{x}$ as
\begin{align}
    \text{mmse}(\alpha) = \lim\limits_{d\to\infty}\frac{1}{d}\mathbb{E}||\hat{\vec{x}}^{\text{opt}}-\vec{x}^{\star}||^{2}_{2} = \rho_{x} - q_x^{\star}
    \label{eq:mmse}
\end{align}
\noindent where $\rho_{x}$ is the second moment of $P_{x}$ and the scalar parameter $q_{x}^{\star}\in[0,\rho_{x}]$ is the solution of the following \emph{free energy} extremisation problem
\begin{align}
   \Phi = -\lim\limits_{d\to\infty}\frac{1}{d}\mathbb{E}_{y}\log\mathcal{Z}(\vec{y}) =   \underset{q_{x},\hat{q}_{x}}{\extr}\left\{\frac{1}{2}\hat{q}_{x}q_{x}-\alpha\Psi_{y}\left(q_{x}\right)-\Psi_{x}(\hat{q}_{x})\right\}\,,\label{eq:freeen}
\end{align}
\noindent with the so-called potentials $(\Psi_{y}, \Psi_{x})$ given by
\begin{align}
\begin{aligned}
\Psi_{y}(t) &= \mathbb{E}_{\xi}\left[\int\dd y~\mathcal{Z}_{y}\left(y;\sqrt{t}\xi, t\right)\log\mathcal{Z}_{y}\left(y;\sqrt{t}\xi, t\right)\right]\, ,\\
\Psi_{x}(r) &= \lim\limits_{d\to\infty}\frac{1}{d} ~\mathbb{E}_{\xi}\left[\mathcal{Z}_{x}(\sqrt{r}\xi,r)\log{\mathcal{Z}_{x}(\sqrt{r}\xi,r)}\right]\, ,
\end{aligned}
\end{align}
\noindent where $\xi\sim\mathcal{N}(0,1)$ and $\mathcal{Z}_{y}$, $\mathcal{Z}_{x}$ are the normalisations of the auxiliary distributions
\begin{align}
Q_{y}\left(x; y,\omega,V\right) &= \frac{1}{\mathcal{Z}_{y}(y;\omega, V)}\frac{e^{-\frac{1}{2V}(x-\omega)^2}}{\sqrt{2\pi V}}\delta\left(y-\varphi(x)\right)\,, \\
Q_{x}\left(\vec{x};B,A\right) &=\frac{P_{x}(\vec{x})}{\mathcal{Z}_{x}(B,A)}\prod\limits_{j=1}^{d}e^{-\frac{A}{2}x_{j}^2+B x_{j}}\,. \nonumber
\end{align}
Note that this expression is valid for arbitrary distribution $P_{x}$, as long as the limit in $\Psi_{x}$ is well-defined. In particular, it reduces to the known result in \cite{krzakala2012probabilistic,Barbier2017c} when $P_{x}$ factorises. In principle, for correlated $P_x$ such as in the generative model of eq.~\eqref{eq:Px} computing $\Psi_{x}$ is itself a hard problem. However, we can see eq.~\eqref{eq:Px} as a chain of generalised linear models. In the limit where $k_{l}\to\infty$ with $\rho = k/d = \Theta(1)$, $L= \Theta(1)$ and $\beta_{l}=k_{l+1}/k_{l} = \Theta(1)$ we can apply the observation above iteratively, layer-wise, up to the input layer for which $P_{z}$ factorises - and is easy to compute. This yields \cite{manoel2017multi}
\begin{multline}
\Phi = \underset{q_{x},\hat{q}_{x}, \{q_{l},\hat{q}_{l}\}}{\extr}\left\{-\frac{1}{2}\hat{q}_{x}q_{x}-\frac{\rho}{2}\sum\limits_{l=1}^{L}\beta_{l}q_{l}\hat{q}_{l}+\alpha\Psi_{y}\left(q_{x}\right) \right. \\ \left. + \rho\sum\limits_{l=2}^{L}\beta_{l}\Psi^{(l)}_{\out}\left(\hat{q}_{l}, q_{l-1}\right)+\Psi^{(L+1)}_{\out}(\hat{q}_{x}, q_{L})+\rho\Psi_{z}\left(\hat{q}_{z}\right)\right\}\label{eq:freeen:multilayer},
\end{multline}
\noindent where we have introduced the additional potentials $(\Psi_{\out}, \Psi_{z})$
\begin{align}
\begin{aligned}
\Psi^{(l)}_{\out}(r,s) &= \mathbb{E}_{\xi,\eta}\left[\mathcal{Z}^{(l)}_{\out}(\sqrt{r}\xi,r,\sqrt{s}\xi, \rho_{l-1}-s)\log{\mathcal{Z}^{(l)}_{\out}(\sqrt{r}\xi,r,\sqrt{s}\xi, \rho_{l-1}-s)}\right]\,, \\
\Psi_{z}(t) &= \mathbb{E}_{\xi}\left[\mathcal{Z}_{z}(\sqrt{t}\xi,t)\log{\mathcal{Z}_{z}(\sqrt{t}\xi,t)}\right]\,,
\end{aligned}
\end{align}
\noindent defined in terms of the following auxiliary distributions
\begin{align}
\begin{aligned}
Q^{(l)}_{\out}(x,z;B,A,\omega,V) &= \frac{e^{-\frac{A}{2}x^2+B x}}{\mathcal{Z}_{\out}(B,A,\omega,V)} \frac{e^{-\frac{1}{2V}\left(z-\omega\right)^2}}{\sqrt{2\pi V}}P^{(l)}_{\out}(x|z)\,,\\
Q_{z}\left(z; B,A\right) &= \frac{e^{-\frac{A}{2}z^2+B z}}{\mathcal{Z}_{z}(B,A)}~P_{z}(z)\,,
\end{aligned}	
\end{align}
\noindent and with $\rho_{l}$ the second moment of the hidden variable $\vec{h}^{(l)}$.

These predictions, that have also been derived with different heuristics in \cite{reeves2017additivity}, were rigorously proven for two-layers in \cite{gabrie2018entropy}, while deeper architectures requires additional assumptions on the concentration of the free energies to be under a rigorous control. Eq.~\eqref{eq:freeen:multilayer} thus reduces the asymptotics of the high-dimensional estimation problem to a low-dimensional extremisation problem over the $2(L+1)$ variables $(q_{x}, \hat{q}_{x}, \{q_{l},\hat{q}_{l}\}_{l=1}^L)$, allowing for a mathematically sound and rigorous investigation. These parameters are also known as the \emph{overlaps}, since they parametrise the overlap between the Bayes-optimal estimator and ground-truth signal at each layer. Solving eq.~\eqref{eq:freeen} provides two important statistical thresholds: the \emph{weak recovery} threshold $\alpha_{c}$ above which better-than-random (i.e. ${\rm{mmse}} < \rho_{x}$) reconstruction becomes theoretically possible and the \emph{perfect reconstruction} threshold, above which perfect signal recovery (i.e. when ${\rm{mmse}} = 0$) becomes possible.

Interestingly, the free energy eq.~\eqref{eq:freeen:multilayer} also provides information about the algorithmic hardness of the problem. The above extremisation problem is closely related the state evolution of the AMP algorithm for this problem, as derived in \cite{manoel2017multi}, and generalized in \cite{fletcher2018inference}. It is conjectured to provide the best polynomial time algorithm for the estimation of $\vec{x}^{\star}$ in our considered setting. Specifically, the algorithm reaches a mean-squared error that corresponds to the local extremiser reached by gradient descent in the function (\ref{eq:freeen:multilayer}) starting with uninformative initial conditions.  

While so far we summarised results that follow from previous works, these results were up to our knowledge not systematically evaluated and analysed for the linear estimation and phase retrieval with generative priors. This analysis and its consequences is the object of the rest of this paper and constitutes the original contributions of this work. 

\subsection{Weak recovery threshold}
Solutions for the extremisation in eq.~\eqref{eq:freeen:multilayer} can be found by solving the fixed point equations, obtained by taking the gradient of eq.~\eqref{eq:freeen:multilayer} with respect of the parameters $(q_{x},\hat{q}_{x}, \{q_{l}, \hat{q}_{l}\}_{l=1}^{L})$:
\begin{align}
\begin{cases}
\hat{q}_{x} =\alpha\Lambda_{y}\left(q_{x}\right)	\\
\hat{q}_{L} =\beta_{L}\Lambda_{\out}\left(\hat{q}_{x}, q_{L}\right)	\\
\hat{q}_{L-1} =\beta_{L-1}\Lambda_{\out}\left(\hat{q}_{L}, q_{L-1}\right)\\
\hspace{1cm}\vdots\\
\hat{q}_{l} =\beta_{l} \Lambda_{\out}\left(\hat{q}_{l+1}, q_{l}\right)\\
\hspace{1cm}\vdots\\
\hat{q}_{z} =\beta_{1} \Lambda_{\out}\left(\hat{q}_{2},q_{z}\right)
\end{cases}&&
\begin{cases}
	q_{x} = \Lambda_{x}\left(\hat{q}_{x}, q_{L}\right)\\
    q_{L} = \Lambda_{x}\left(\hat{q}_{L}, q_{L-1}\right)\\
\hspace{1cm}\vdots\\
{q}_{l} = \Lambda_{x}\left(\hat{q}_{l},q_{l-1}\right) \\
\hspace{1cm}\vdots\\
q_{z} =\Lambda_{z}\left(\hat{q}_{z}\right)
\end{cases}
	\label{eq:SE_uu}
\end{align}
\noindent where $\Lambda_{y}(t) = 2~\partial_{t}\Psi_{y}(t)$, $\Lambda_{z}(t) = 2~\partial_{t}\Psi_{z}(t)$, $\Lambda_{x}(t) = 2~\partial_{r}\Psi_{\out}(r,s)$, $\Lambda_{\out}(t) = 2~\partial_{s}\Psi_{\out}(r,s)$.

The weak recovery threshold $\alpha_{c}$ is defined as the value above which one can estimate $\vec{x}^{\star}$ better than a random draw from the prior $P_{x}$. In terms of the mmse it is defined as
\begin{align}
\alpha_{c} = \underset{\alpha\geq 0}{\rm{argmax}	}\{{\rm{mmse}}(\alpha)=\rho_{x}\}.
\end{align}
From eq.~\eqref{eq:mmse}, it is clear that an uninformative solution ${\rm{mmse}}=\rho_{x}$ of eq.~\eqref{eq:freeen:multilayer} corresponds to a fixed point $q_{x} = 0$. For both the phase retrieval and linear estimation, evaluating the right-hand side of eqs.~\eqref{eq:SE_uu} at $q_{x}=0$ we can see that $\hat{q}^{\star}_{x}=0$ is a fixed point if $\sigma$ is an odd function and if
\begin{align}
\mathbb{E}_{P_{z}} \[ z \] = 0, \hspace{ 1.5cm }\textrm{ and }\hspace{ 1.5cm } \mathbb{E}_{Q^{(l), 0}_{\out}}\[ x \] = 0 \,,
\label{eq:stability_conditions} 	
\end{align}
\noindent where $Q^{(l), 0}_{\out}(x,z) = Q^{(l)}_{\out}(x,z;0,0,0,\rho_{l-1})$. These conditions reflect the intuition that if the prior~$P_{z}$ or the likelihoods $P_{\out}^{(l)}$ are biased towards certain values, this knowledge helps the statistician estimating better than a random guess. If these conditions are satisfied, then $\alpha_{c}$ can be obtained as the point for which the fixed point $q_{x}=0$ becomes unstable. The stability condition is determined by the eigenvalues of the Jacobian of eqs.~\eqref{eq:SE_uu} around the fixed point $(q^{\star}_{x}, \hat{q}^{\star}_{x}, \{q^{\star}_{l}, \hat{q}^{\star}_{l}\}_{l=1}^L) = 0$. More precisely, the fixed point becomes unstable as soon as one eigenvalue of the Jacobian is bigger than one. Expanding the update functions around the fixed point and using the conditions in eq.~\eqref{eq:stability_conditions},
\begin{align}
\Lambda_{y}(t) &\underset{t\ll 1}{=} \frac{1}{\rho^2_{x}}\int\dd y~\mathcal{Z}_{y}(y;0,\rho_{x})\left(\mathbb{E}_{Q^{0}_{y}}[\rho_{x}-x^2]\right)^2 t+\mO \left(t^{3/2}\right) \notag \,,\\
\Lambda^{(l)}_{x}(r,s) &\underset{r,s\ll 1}{=} \left(\mathbb{E}_{Q^{(l), 0}_{\out}}[x^2]\right)^2~r+\frac{1}{\rho_{l-1}^2}\left(\mathbb{E}_{Q^{(l), 0}_{\out}}[xz]\right)^2 s+\mO \left(r^{3/2},s^{3/2}\right)\,, \\
\Lambda^{(l)}_{\out}(r,s) &\underset{r,s\ll 1}{=}\left(\mathbb{E}_{Q^{(l), 0}_{\out}}[xz]\right)^2  r+\frac{1}{\rho_{l-1}^2}	\left(\mathbb{E}_{Q^{(l), 0}_{\out}}[z^2]-\rho_{l-1}\right)^2 s+\mO\left(r^{3/2},s^{3/2}\right)\notag \,, \\
\Lambda_{z}(t) &\underset{t\ll 1}{=}\left(\mathbb{E}_{P_z}[z^2]\right)^2 t+\mO\left(t^{3/2}\right)\,. \notag
\end{align}
For a generative prior with depth $L$, the Jacobian is a cumbersome sparse $(L+1) \times (L+1)$ matrix, with all the entries given by the six partial derivatives above. For the sake of conciseness we only write it here for $L=1$:
\begin{align}
\begin{pmatrix}
0 & \left(\mathbb{E}_{Q^{0}_{\out}}\[x^2\]\right)^2 &\frac{1}{\rho_{z}^{2}}\left(\mathbb{E}_{Q^{0}_{\out}}\[xz\] \right)^2 & 0\\ 	
\frac{\alpha}{\rho_{x}^2}\int\dd y~\mathcal{Z}^{0}_{y} \left(\mathbb{E}_{Q^{0}_{y}}\[\rho_{x}-x^2\]\right)^2 & 0 & 0 & 0\\
0 & 0 & 0 & \left(\mathbb{E}_{P_{z}}\[z^2\]\right)^2\\
0 & \beta\left(\mathbb{E}_{Q^{0}_{\out}}\[xz\]\right)^2 &\frac{\beta}{\rho_{z}^2}\left(\mathbb{E}_{Q^{0}_{\out}}\[z^2\]-\rho_{z}\right)^2 & 0
\end{pmatrix}.
\end{align}
Note that this holds for any choice of $P_{\out}^{(l)}$ and latent space distribution $P_{z}$, as long as conditions eq.~\eqref{eq:stability_conditions} hold. For the phase retrieval with a linear generative model for instance $P^{(l)}(x|z) = \delta(x-z)$, we find $\alpha_{c} = \frac{1}{2}\frac{1}{1+\rho^{-1}}$. For a linear network of depth $L$ this generalises to
\begin{align}
\alpha_{c} = \frac{1}{2}\left(1+\sum\limits_{l=1}^{L}\prod\limits_{k=0}^{l-1}\beta_{L-k}\right)^{-1}.
\label{eq:weak_recovery}
\end{align}
The linear estimation problem has exactly the same threshold, but without the global $1/2$ factor. Since $\rho, \beta_{l}\geq 0$, it is clear that $\alpha_{c}$ is decreasing in the depth $L$ of the network. This analytical formula is verified by numerically solving eqs.~\eqref{eq:SE_uu}, see Figs.~\ref{main:phase_diagramm_CS} and \ref{main:phase_diagramm_PR}. For other choices of activation satisfying condition \eqref{eq:stability_conditions} (e.g. the sign function), we always find that depth helps in the weak recovery of the signal.

\subsection{Perfect recovery threshold}
We now turn our attention to the perfect recovery threshold, above which perfect signal reconstruction becomes statistically possible. Formally, it can be defined as
\begin{align}
\alpha_{\IT} &= \underset{\alpha\geq 0}{\text{argmin}}\{\text{mmse}(\alpha) = 0\},
\end{align}
\noindent and corresponds to the global minimum of the free energy in eq.~\eqref{eq:freeen:multilayer}. Numerically, the perfect recovery threshold is found by solving the fixed point equations \eqref{eq:SE_uu} from an informed initialisation $q_{x}\approx \rho_{x}$, corresponding to $\text{mmse} \approx 0$ according to eq.~\eqref{eq:mmse}. The resulting fixed point is then checked to be a minimiser of the free energy eq.~\eqref{eq:freeen:multilayer}. Different from $\alpha_{c}$, it cannot be computed analytically for an arbitrary architecture. However, for the compressed sensing and phase retrieval problems with $\sigma\in\{\text{linear}, \text{ReLU}\}$ generative priors, $\alpha_{\IT}$ can be analytically computed by generalising a simple argument based on the invertibility of the linear system of equations at each layer, originally used in the usual compressive sensing \cite{candes2006near, tao2009}. 

First, consider the linear estimation problem with a deep linear generative prior, i.e. $\vec{y} = \mat{A} \vec{x}^{\star}\in\mathbb{R}^{n}$ with $\vec{x}^{\star} = \mat{W}^{(L)}\dots \mat{W}^{(1)}\vec{z}\in\mathbb{R}^{d}$ and $\mat{A}, \{\mat{W}^{(l)}\}_{l=1}^L$ i.i.d. Gaussian matrices, that are full rank with high probability. For $n > d$, the system $\vec{y} = \mat{A} \vec{x}^{\star}$ is overdetermined as there are more equations than unknowns. Hence the information theoretical threshold has to verify $\alpha_{\IT} = \frac{n_{\IT}}{d} \leq 1$. For $L=0$ (i.e. $\vec{x}^{\star}$ is Gaussian i.i.d.), we have exactly $\alpha_{\IT}^{(0)}=1$ as the prior does not give any additional information for solving the linear system. 
For $L\geq 1$ though, at each level $l \in [1:L]$, we need to solve successively $\vec{h}^{(l)} \in \bbR^{k_l}$ in the linear system $\vec{y} = \mat{A} \mat{W}^{(L)} \cdots \mat{W}^{(l)} \vec{h}^{(l)}$. Again as $\mat{A} \mat{W}^{(L)} \cdots \mat{W}^{(l)} \in\bbR^{ n \times k_l}$, if $n > k_l$ the system is over-constrained. Hence the information theoretical threshold for this equation is such that $\forall l \in [1:L],  n_{\IT}^{(l)} \leq k_l  \Leftrightarrow \alpha_{\IT}^{(l)} \leq \prod\limits_{k=1}^{l} \frac{1}{\beta_{L-k+1}} $. And note that $\rho \equiv \prod\limits_{k=1}^{L} \frac{1}{\beta_{L-k+1}}$.
Hence, the information theoretical threshold is obtained by taking the smallest of the above values $\alpha_{\IT}^{(l)}$:
			\begin{align}
			 	\alpha_{\IT} = \min_{l\in [0:L]} \alpha_{\IT}^{(l)} = \min \left(1, \left\{ \prod\limits_{k=1}^{l} \frac{1}{\beta_{L-k+1}}  \right\}_{l=1}^{L-1}, \rho \right)\,.
			 	\label{eq:alphaIT}
			\end{align}

This result generalises to the real-valued phase retrieval problem. First, we note that by the data processing inequality taking $\vec{y} = |\mat{A}\vec{x}^{\star}|$ cannot increase the information about $\vec{x}^{\star}$, and therefore the transition in phase retrieval cannot be {\it better} than for compressed sensing. Secondly, an inefficient algorithm exists that achieve the same performance as compressed sensing for the real valued phase retrieval: one just needs to try all the possible $2^m$ assignments for the sign, and then solve the corresponding compressed sensing problem. This strategy that will work as soon as the compressed sensing problem is solvable. Eq.~(\ref{eq:alphaIT}) is thus valid for the real phase retrieval problem as well.

One can finally generalise this analysis for a non-linear generative prior with ReLU activation at each layer, i.e. $\vec{x}^{\star} = \text{relu}\left(\mat{W}^{(L)}\text{relu}\left(\cdots \mat{W}^{(1)}\vec{z}\right)\cdots\right)$. Noting that on average $\vec{x}$ has half of zero entries and half of {\iid} Gaussian entries, the system can be reorganised and simplified $\vec{y} = \td{\mat{A}} \td{\vec{x}}$, with $\td{\vec{x}}\in \bbR^{d/2}$ the extracted vector of $\vec{x}$ with on average $d/2$ strictly positive entries and the corresponding reduced matrix $\td{\mat{A}} \in \bbR^{n \times d/2}$, is over-constrained for $n > d/2 $ and hence the information theoretical threshold verifies $\alpha_{\IT} = \frac{n_{\IT}}{d} \leq \frac 12 $. 
Noting that this observation remains valid for generative layers, we will have on average at each layer an input vector $\vec{h}^{(l)}$ with half of zero entries and half of Gaussian distributed entries - except at the very first layer for which the input $\vec{z}\in\mathbb{R}^{k}$ is dense. Repeating the above arguments yields the following perfect recovery threshold
			\begin{align}
			 	\alpha_{\IT} =  \text{min}\left(\frac{1}{2}, \left\{ \frac{1}{2} \prod\limits_{k=1}^{l} \frac{1}{\beta_{L-k+1}}  \right\}_{l=1}^{L-1}, \rho \right)\,.
			 	\label{eq:alphaIT2}
			\end{align}
	for both the linear estimation and phase retrieval problems. Both these results are consistent  with the solution of the saddle-point eqs.~\eqref{eq:SE_uu} with a informed initialisation, see Figs.~\ref{main:phase_diagramm_PR}-\ref{fig:multilayer:compression}.

\subsection{Algorithmic threshold}
\label{sec:algo}
The discussion so far focused on the statistical limitations for signal recovery, regardless of the cost of the reconstruction procedure. In practice, however, one is concerned with the algorithmic costs for reconstruction. In the high-dimensional regime we are interested, where the number of observations scale with the number of parameters in the model, only (low)-polynomial time algorithms are manageable in practice. Remarkably, the formula in eq.~\eqref{eq:freeen:multilayer} also provides useful information about the algorithmic hindrances for the inverse problems under consideration. Indeed, with a corresponding choice of iteration schedule and initialisation, the fixed point equations eq.~\eqref{eq:freeen:multilayer} are identical to the state evolution describing the asymptotic performance of an associated \emph{Approximate Message Passing} (AMP) algorithm \cite{manoel2017multi,fletcher2018inference}. Moreover, the AMP aforementioned is the \emph{best known} polynomial time algorithm for the estimation problem under consideration, and it is conjectured to be the optimal polynomial algorithm in this setting.

\begin{figure}[htb!]
	\centering
		\includegraphics[scale=0.45]{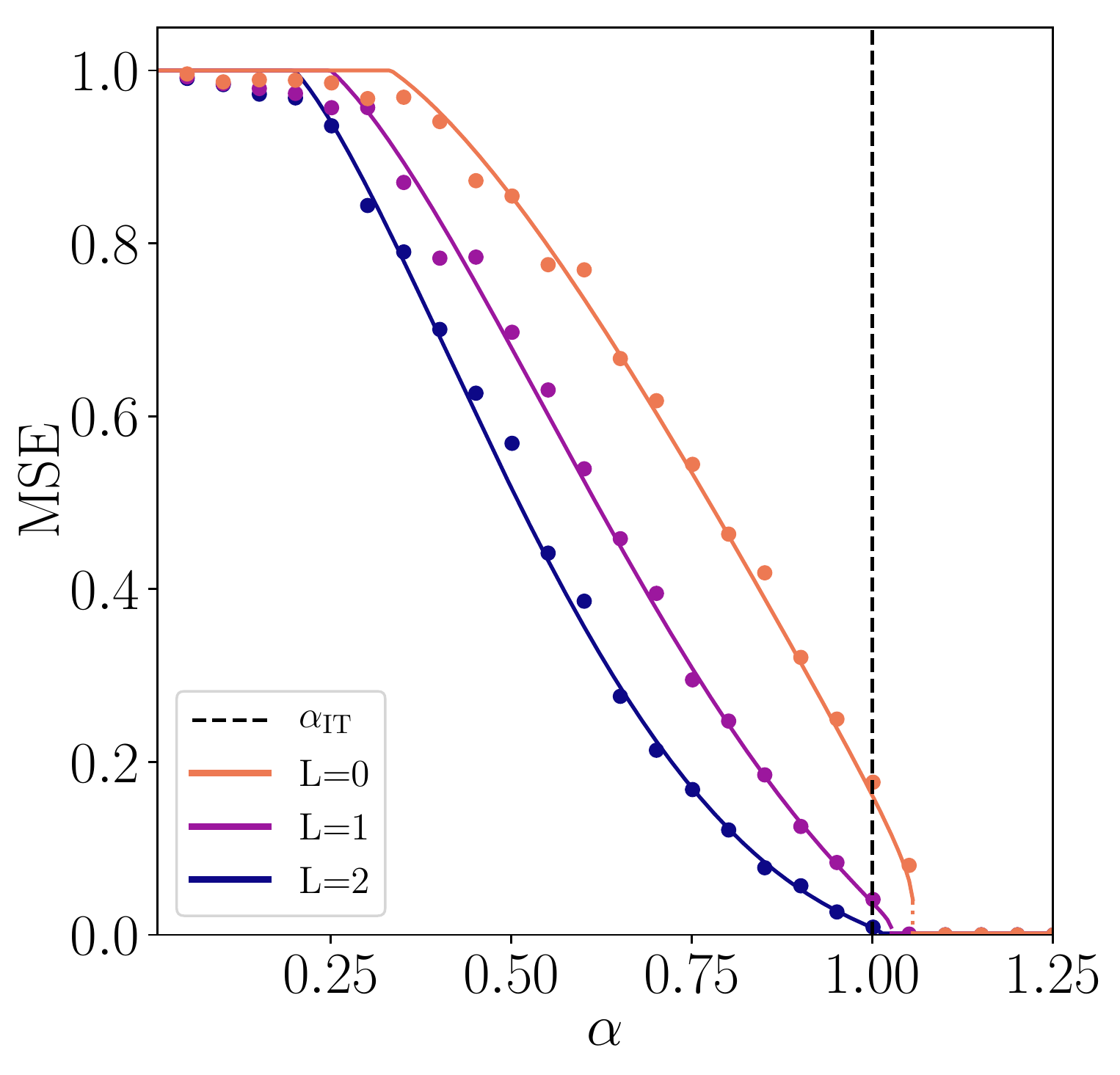}%
		\includegraphics[scale=0.45]{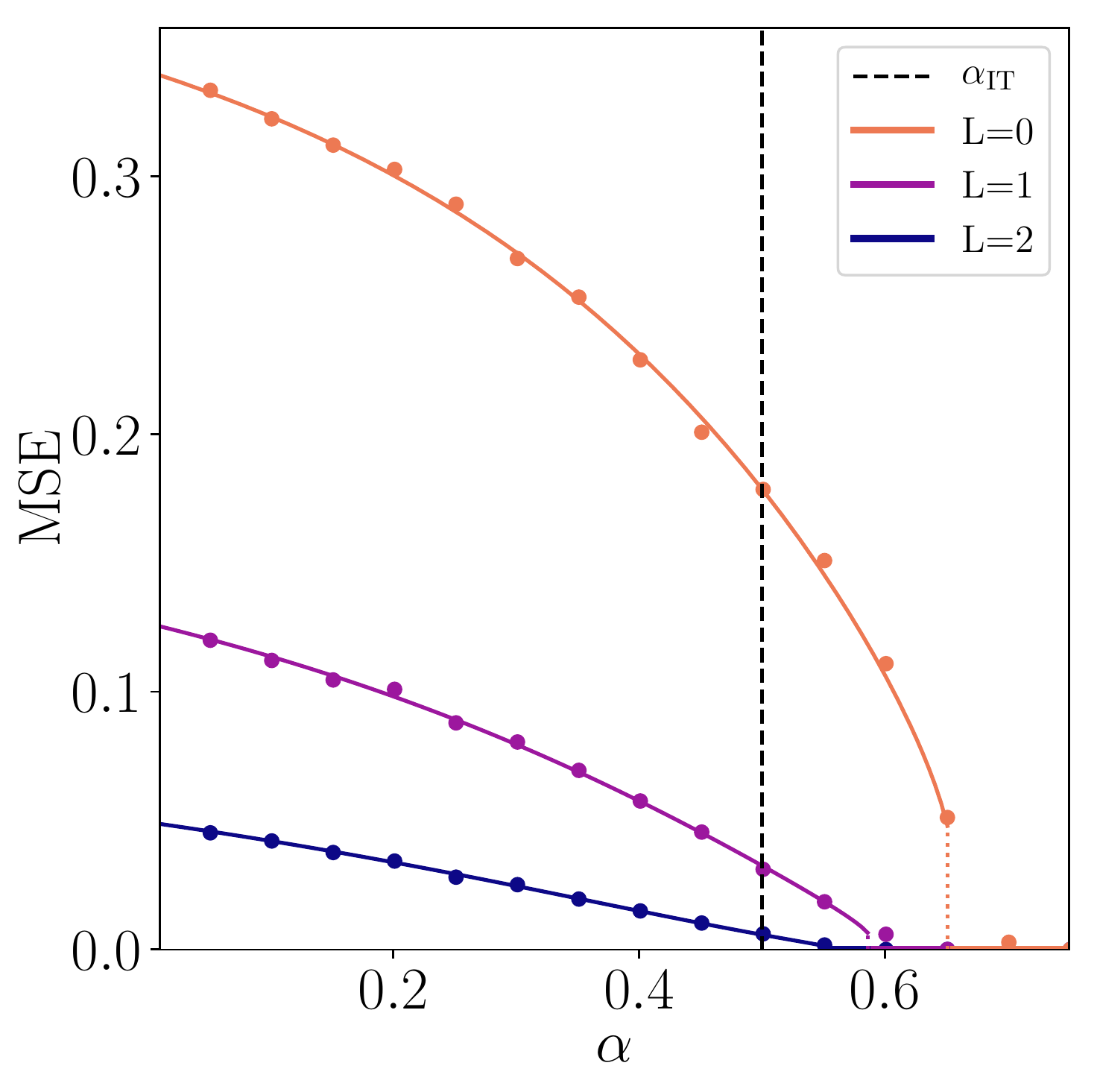}
	\caption{Mean squared error obtained by running the AMP algorithm (dots) from \cite{manoel2017multi} and implemented with the \texttt{tramp} package \cite{baker2020tramp}, for $d=2.10^3$ averaged on $10$ samples, compared to the MSE obtained from the state evolution eqs.~\eqref{eq:SE_uu} with uninformative initialisation $q_{x}=q_{l}\approx 0$ (solid line) for the phase retrieval problem with linear (\textbf{left}) and relu (\textbf{right}) generative prior networks. Different curves correspond to different depths $L$, with fixed $\rho = 2$ and layer-wise aspect ratios $\beta_{l} = 1$. The dashed vertical line corresponds to $\alpha_{\IT}$. To illustrate for instance in the linear case (\textbf{left}), $(\alpha_c^{L=0}, \alpha_c^{L=1},\alpha_c^{L=2}) = (1/3, 1/4, 1/5)$, $\alpha_{\IT}=1$ and  $(\alpha_{\alg}^{L=0}, \alpha_{\alg}^{L=1}, \alpha_{\alg}^{L=2}) = (1.056, 1.026 , 1.011)$.
	\label{fig:mses}
	\vspace{-0.2cm}
	}
\end{figure}

The AMP state evolution corresponds to initialising the overlap parameters $(q_{x},q_{l}) \approx 0$ and updating, at each time step $t$ the hat variables $\hat{q}_{x}^{t} = \alpha\Lambda_{y}(q_{x}^{t})$ before the overlaps $q_{x}^{t+1} = \Lambda_{x}(\hat{q}_{x}^{t}, q_{L}^{t})$, etc. In Fig.~\ref{fig:mses} we illustrate this equivalence by comparing the MSE obtained by iterating eqs.~\eqref{eq:SE_uu} with the averaged MSE obtained by actually running the AMP algorithm from \cite{manoel2017multi} for a specific architecture and implemented with the \texttt{tramp} python package \cite{baker2020tramp}. In particular even though the AMP state evolution is not yet rigorously proven, we see a strong agreement of our analytical results with AMP simulations.

Note that, by construction, the performance of the Bayes-optimal estimator corresponds to the global minimum of the scalar potential in eq.~\eqref{eq:freeen:multilayer}. If this potential is convex, eqs.~\eqref{eq:SE_uu} will converge to the global minimum, and the asymptotic performance of the associated AMP algorithm will be optimal. However, if the potential has also a local minimum, initialising the fixed point equations will converge to the different minima depending on the initialisation. In this case, the MSE associated to the AMP algorithm (corresponding to the local minimum) differs from the Bayes-optimal one (by construction the global minimum). In the later setting, we define the \emph{algorithmic threshold} as the threshold above which AMP is able to perfectly reconstruct the signal - or equivalently for which $\text{mmse} = 0$ when eqs.~\eqref{eq:SE_uu} are iterated from $q_{x}^{t=0}=q_{l}^{t=0}=\epsilon\ll 1$. Note that by definition $\alpha_{\IT} < \alpha_{\alg}$, and we refer to $\Delta_{\alg} = \alpha_{\alg} - \alpha_{\IT}$ as the algorithmic gap. See Fig.~\ref{fig:landscape} for an illustration of the evolution of the free energy landscape for increasing $\alpha$.

Studying the existence of an algorithmic gap for the linear estimation and phase retrieval problems, and how it depends on the architecture and depth of the generative prior, is the subject of the next section.

%% file: files/phase_diagrams.tex
\section{Phase diagrams}
\label{sec:phase_diagrams} 
In this section we summarise the previous discussions in plots in the $(\rho,\alpha)$-plane, hereafter named \emph{phase diagrams}. Phase diagrams quantify the quality of signal reconstruction for a fixed architecture $(\beta_1, \dots, \beta_{L-1})$ \footnote{Note that $\beta_{L}$ is fixed from the knowledge of $(\rho, \beta_1, \dots, \beta_{L-1})$.} as a function of the compression $\rho$. Moreover, it allows a direct visual comparison between the phase diagram for a sparse Gaussian prior and the multi-layer generative prior. For both the phase retrieval and compressed sensing problems we distinguish the following regions of parameters limited by the thresholds of sec.~\ref{sec:it_analysis}: 
\begin{itemize}
    \item 
    {\it Undetectable} region where the best achievable error is as bad as a random guess from the prior as if no measurement $\vec{y}$ were available. Corresponds to $\alpha < \alpha_{c}$.
    \item 
    {\it Weak recovery} region where the optimal reconstruction error is better than the one of a random guess from the prior, but exact reconstruction cannot be achieved. Corresponds to $\alpha_c < \alpha < \alpha_{\IT}$.
    \item 
    {\it Hard} region where exact reconstruction can be achieved information-theoretically, but no efficient algorithm achieving it is known. Corresponds to $\alpha_{\IT}<\alpha<\alpha_{\alg}$ 
   \item 
    The so-called {\it easy} region where the aforementioned AMP algorithm for this problem achieves exact reconstruction of the signal. Corresponds to $\alpha > \alpha_{\alg}$.
\end{itemize}
As already explained, we locate the corresponding phase transitions in the following manner:
for the weak recovery threshold $\alpha_c$, we notice that the fixed point corresponding to an error as bad as a random guess corresponds to the values of the order parameters $q_{x}, q_{l} = 0$. This is an extremiser of the free energy (\ref{eq:freeen}) when the prior $P_{z}$ has zero mean and the non-linearity $\varphi$ is an even function. This condition is satisfied for both the linear estimation and the phase retrieval problem with linear generative priors that leads to zero-mean distributions on the components of the signal, but is not achieved for a generative prior with ReLU activation, since it biases estimation.
In case this uninformative fixed point exists, we investigate its stability under the state evolution of the AMP algorithm, thus defining the threshold $\alpha_c$. For $\alpha<\alpha_c$ the fixed point is stable, implying the algorithm is not able to find an estimator better than random guess. In contrast, for $\alpha>\alpha_c$ the AMP algorithm provides an estimator better than random guess. For phase retrieval with linear generative model in the setting of the present paper, this analysis leads to the threshold derived in eq.~\eqref{eq:weak_recovery}.
If there exists a region where the performance of the AMP algorithm and the information-theoretic one do not agree we call it the {\it hard} region. The hard region is delimited by threshold $\alpha_{\IT}$ and $\alpha_{\alg}$. 

\begin{figure}[t]
	\centering
	\begin{tikzpicture}[scale=0.20]
    \tikzstyle{dot}=[circle,minimum size=10pt, scale=0.4]
    \tikzstyle{green}=[fill=green!70!black]
    \tikzstyle{orange}=[fill=orange!90!black]
    \tikzstyle{yellow}=[fill=yellow!95!black]
    \tikzstyle{red}=[fill=red!90!black]
    \tikzstyle{annot}=[text width=3cm, text centered, font=\footnotesize]
    \tikzstyle{line}=[-,thick, black]
	\draw [line] (0,0) to [out=-65,in=180] (3.5,-2) to [out=0,in=200] (9,1) ;
	\node[dot, green] at (3.5,-2) {}; 
	\draw [line] (11,0) to [out=-45,in=180] (13.5,-2) to [out=0,in=-180] (15.5,0) to [out=0,in=-180] (16,0) to [out=10,in=220] (19,2)  ;
	\node[dot, green] at (13.5,-2) {}; 
	\node[dot, orange] at (16,0) {}; 
	\node[annot, below = 0.25cm] at (15, -2) {$\alpha_{\rm sp}$}  ;
	\draw [line] (15, -6.5) to (15, -5.5) ;
	\draw [line] (21,0) to [out=-45,in=180] (23.5,-2) to [out=0,in=180] (25,0) to [out=0,in=-180] (26.5,-1) to [out=0,in=180] (29,0);
	\node[dot, green] (BO) at (23.5,-2) {}; 
	\node[dot, orange] (BO) at (26.5,-1) {}; 
	\draw[densely dotted] (33.5,-2) to (36.5,-2);
	\draw [line] (31,0) to [out=-45,in=180] (33.5,-2) to [out=0,in=180] (35,0) to [out=0,in=-180] (36.5,-2) to [out=0,in=180] (39,0);
	\node[dot, red]  at (33.5,-2) {}; 
	\node[dot, green] at (36.5,-2) {}; 
	\node[annot, below = 0.25cm] at (35, -2) {$\alpha_{\rm IT}$}  ;
	\draw [line] (35, -6.5) to (35, -5.5) ;
	\draw [line] (41,0) to [out=-45,in=180] (43.5,-1) to [out=0,in=180] (45,0) to [out=0,in=-180] (46.5,-2) to [out=0,in=180] (49,0);
	\node[dot, red] at (43.5,-1) {}; 
	\node[dot, green] at (46.5,-2) {}; 
	\draw [line] (51,1) to [out=-70,in=180] (53.5,0) to [out=0,in=-180] (56.5,-2) to [out=0,in=180] (59,0);
	\node[dot, red] at (53, 0) {}; 
	\node[dot, green] at (56.5,-2) {}; 
	\node[annot, below = 0.25cm] at (55, -2) {$\alpha_{\rm alg}$} ;
	\draw [line] (55, -6.5) to (55, -5.5) ;
	\draw [line] (61,1) to [out=0,in=-180] (66.5,-2) to [out=0,in=180] (69,0);
	\node[dot, green] at (66.5,-2) {};
	\draw[draw=none, orange] (0,-6) to (0,-6.5) to (35,-6.5) to (35, -6) to (0, -6);
	\draw[draw=none, yellow] (35,-6) to (35,-6.5) to (55,-6.5) to (55, -6) to (35, -6);
	\draw[draw=none, green] (55,-6) to (55,-6.5) to (70,-6.5) to (70, -6) to (55, -6);
	\draw [->,thick, black] (0,-6) to (70,-6);
	\node[text centered, font=\footnotesize, below = 0.25cm] at (17.5, -6.2) {\textbf{Weak recovery}}  ;
	\node[annot, below = 0.25cm] at (45, -6.2) {\textbf{Hard}}  ;
	\node[annot, below = 0.25cm] at (62.5, -6.2) {\textbf{Easy}}  ;
	\end{tikzpicture}
\caption{Illustration of the free energy landscape as a function of the overlap with the ground truth solution, when one increases $\alpha$. For small $\alpha<\alpha_{\rm sp}$, there exists a unique global minimum, whose overlap with the solution is small  (high MSE). 
	At $\alpha=\alpha_{\rm sp}$, a \emph{local} minimum (orange dot) with higher overlap (small MSE) appears. By definition, the global minimum corresponds to the MMSE of the problem, which is the MSE attained by the Bayes-optimal estimator (green dot). For $\alpha<\alpha_{\rm IT}$ the accessible solution, i.e the global minimum (green dot) has a high MSE while a better solution exists but has a higher free energy (weak recovery phase). At $\alpha=\alpha_{\rm IT}$ the two minima are global and have the same free energy.  Between $\alpha_{\rm IT} < \alpha <\alpha_{\rm alg}$ (hard phase), the local minimum with higher MSE corresponds to the performance of the AMP estimator (red dot). Above $\alpha_{\rm alg}$ only the small MSE minima survive and the AMP estimator is able to achieve the Bayes-optimal performance (easy phase). 
}
\label{fig:landscape}
\end{figure}
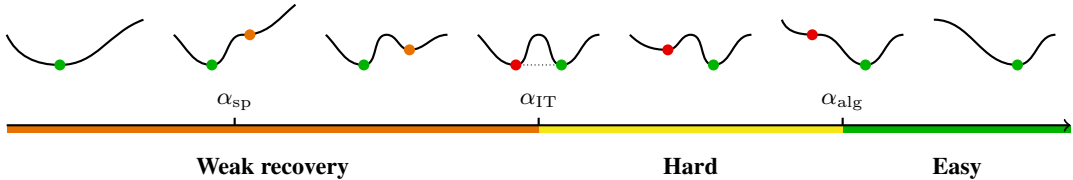

The statistical and algorithmic thresholds defined above admit an alternative and instructive description in terms of free energy landscape, see Fig. ~\ref{fig:landscape}. Consider a fixed $\rho$: for small $\alpha$ the free energy eq.~\eqref{eq:freeen:multilayer} has a single global minimum with small overlap (high MSE) with the ground truth solution $\vec{x}^\star$, referred as the \emph{uninformative} fixed point. 
At a value $\alpha_{\rm sp}$, known as the \emph{first spinodal transition}, a second local minimum appears with higher overlap (smaller MSE) with the ground truth, referred as \emph{informative} fixed point. The later fixed point becomes a global minimum of the free energy at $\alpha_{\IT} > \alpha_{\rm sp}$, while the uninformative fixed point becomes a local minimum. A second spinodal transition occurs at $\alpha_{\alg}$ when the informed fixed point becomes unstable.  
Numerically, the informed and uninformative fixed points can be reached by iterating the saddle-point equations from different initial conditions. When the two are present, the informed fixed point can be reached by iterating from $q_{x}\approx \rho_{x}$, which corresponds to a minimum overlap with the ground truth $\bx^{\star}$, and the uninformative fixed point from $q_{x}\approx 0$, corresponding to no initial overlap with the signal. In the noiseless linear estimation and phase retrieval studied here we observe $\alpha_{\IT} = \alpha_{\rm sp}$. 

\subsection{Single-layer generative prior}
\label{sec:phase_diagrams:single_layer}
First, we consider the case where the signal is generated from a single-layer generative prior, $\vec{x}^{\star} = \sigma(\mat{W}\vec{z})$ with $\vec{z}\sim\mathcal{N}(0,\mat{I}_{k})$. We analyse both the compressed sensing and the phase retrieval problem, for $\sigma\in\{\text{linear}, \text{ReLU}\}$. In this case the only free parameters of the model are $(\rho, \alpha)$, and therefore the phase diagram fully characterises the recovery in these inverse problems. The aim is to compare with the phase diagram of a sparse prior with density $\rho_{s} = \rho$ of nonzero components.

\begin{figure}[htb!]
	\centering
		\includegraphics[scale=0.46]{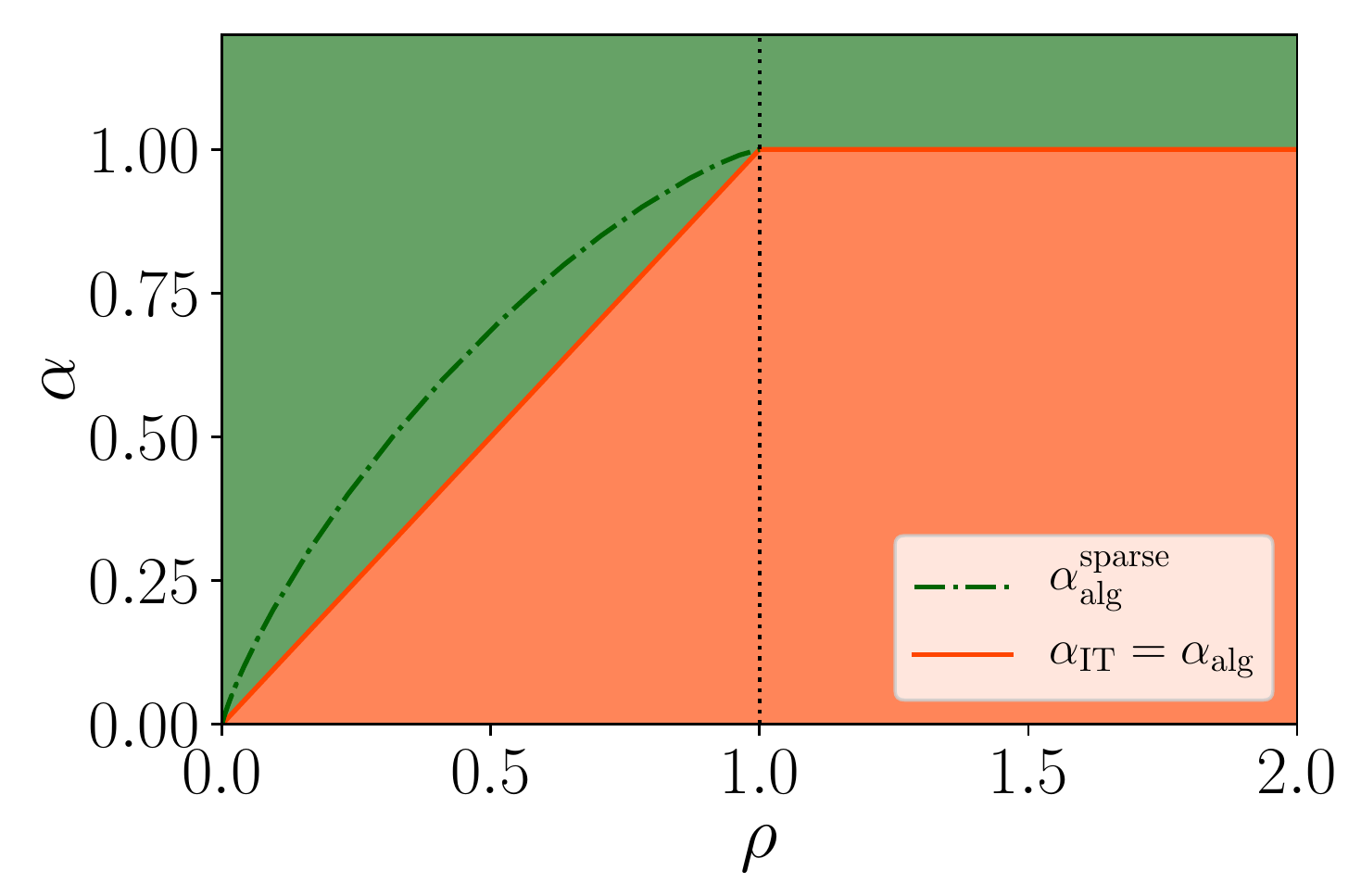}%
		\includegraphics[scale=0.46]{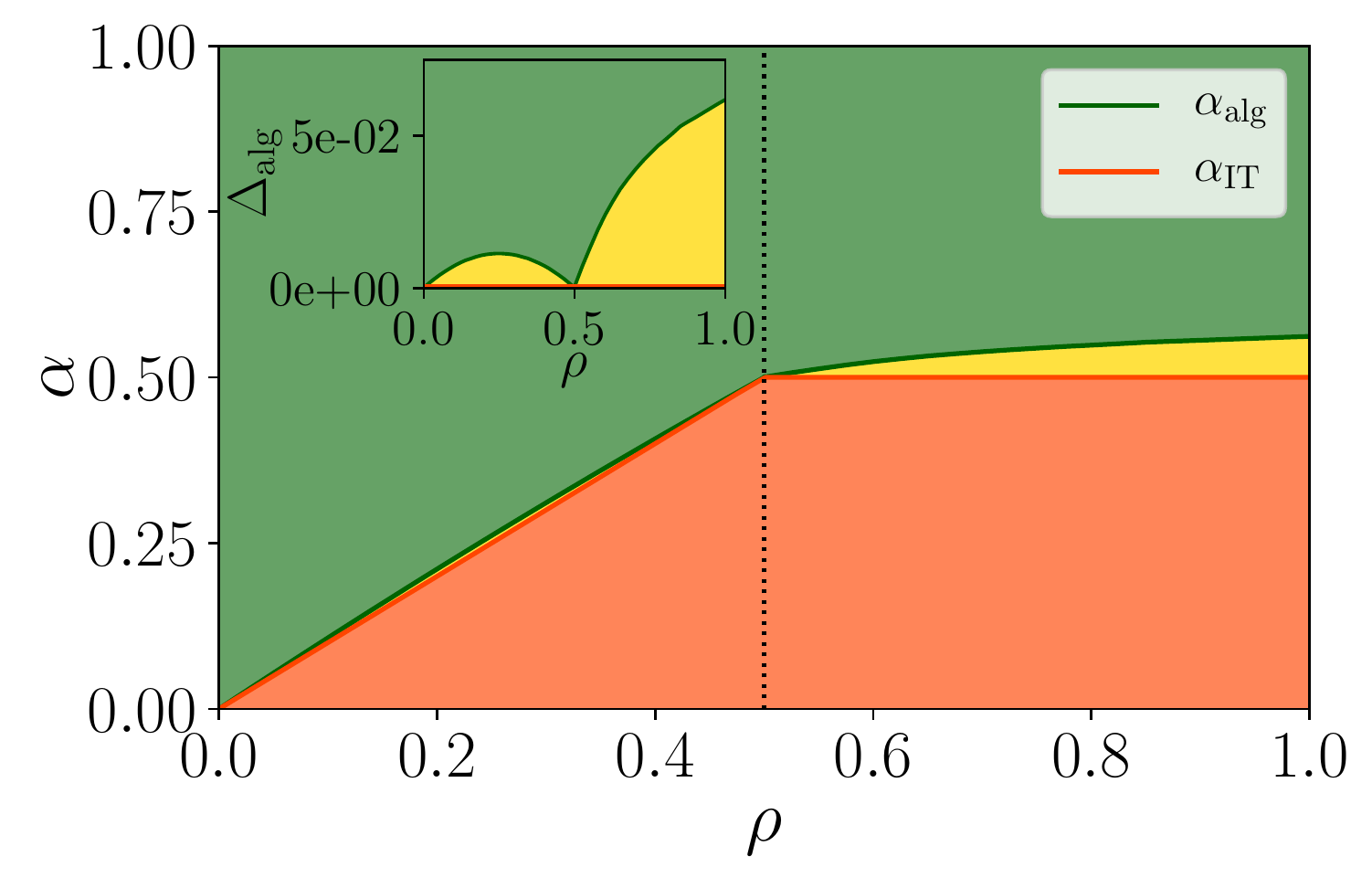}
	\caption{Phase diagrams for the compressed sensing problem with \textbf{(\textbf{left})} linear generative prior and \textbf{(\textbf{right})} ReLU generative prior, in the plane $(\rho,\alpha)$. The $\alpha_{\IT}$ (red line) represents the information theoretic transition for perfect reconstruction and $\alpha_{\rm alg}$ (green line)  the algorithmic transition to perfect reconstruction. In the left part we depict for comparison the algorithmic phase transition for sparse separable prior $\alpha_{\rm alg}^{\rm sparse}$ (dashed-dotted green line). The inset in the right part depicts the difference $\Delta_{\rm alg} = \alpha_{\rm alg} - \alpha_{\IT}$. Colored areas correspond respectively to the \textit{weak recovery} (orange), \textit{hard} (yellow) and \textit{easy} (green) phases. The behaviour of the free energy landscape for increasing $\alpha$ and fixed $\rho$ is illustrated in Fig.~\ref{fig:landscape}. 
	\label{main:phase_diagramm_CS}
	\vspace{-0.2cm}
	}
\end{figure}

Fig.~\ref{main:phase_diagramm_CS} depicts the compressed sensing problem with linear (left) and ReLU (right) generative priors. We depict the phase transitions defined above. On the left hand side we compare to the algorithmic phase transition known from \cite{krzakala_statistical-physics-based_2012} for sparse separable prior with fraction $1-\rho$ of zero entries and $\rho$ of Gaussian entries of zero mean presenting an algorithmically hard phase for $\rho < \alpha < \alpha_{\rm alg}^{\rm sparse}(\rho)$. 

In the case of compressed sensing with linear generative prior we do not observe any hard phase and exact recovery is possible for $\alpha \ge \min (\rho,1)$ due to invertibility (or the lack of there-of) of the matrix product~$\mat{A}\mat{W}$. With ReLU generative prior we have $\alpha_{\IT} = \min(\rho,1/2)$ and the hard phase exists and has interesting properties: The $\rho \to \infty$ limit corresponds to the separable prior, and thus in this limit $\alpha_{\rm alg}(\rho \to \infty) = \alpha_{\rm alg}^{\rm sparse}(\rho_s=1/2)$. Curiously we observe $\alpha_{\rm alg} > \alpha_{\IT}$ for all $\rho \in (0,\infty)$ except at $\rho=1/2$. Moreover the size of the hard phase is very small for $\rho < 1/2$ when compared to the one for compressed sensing with separable priors, suggesting that exploring structure in terms of generative models might be algorithmically advantageous over sparsity.

\begin{figure}[htb!]
	\centering
		\includegraphics[scale=0.45]{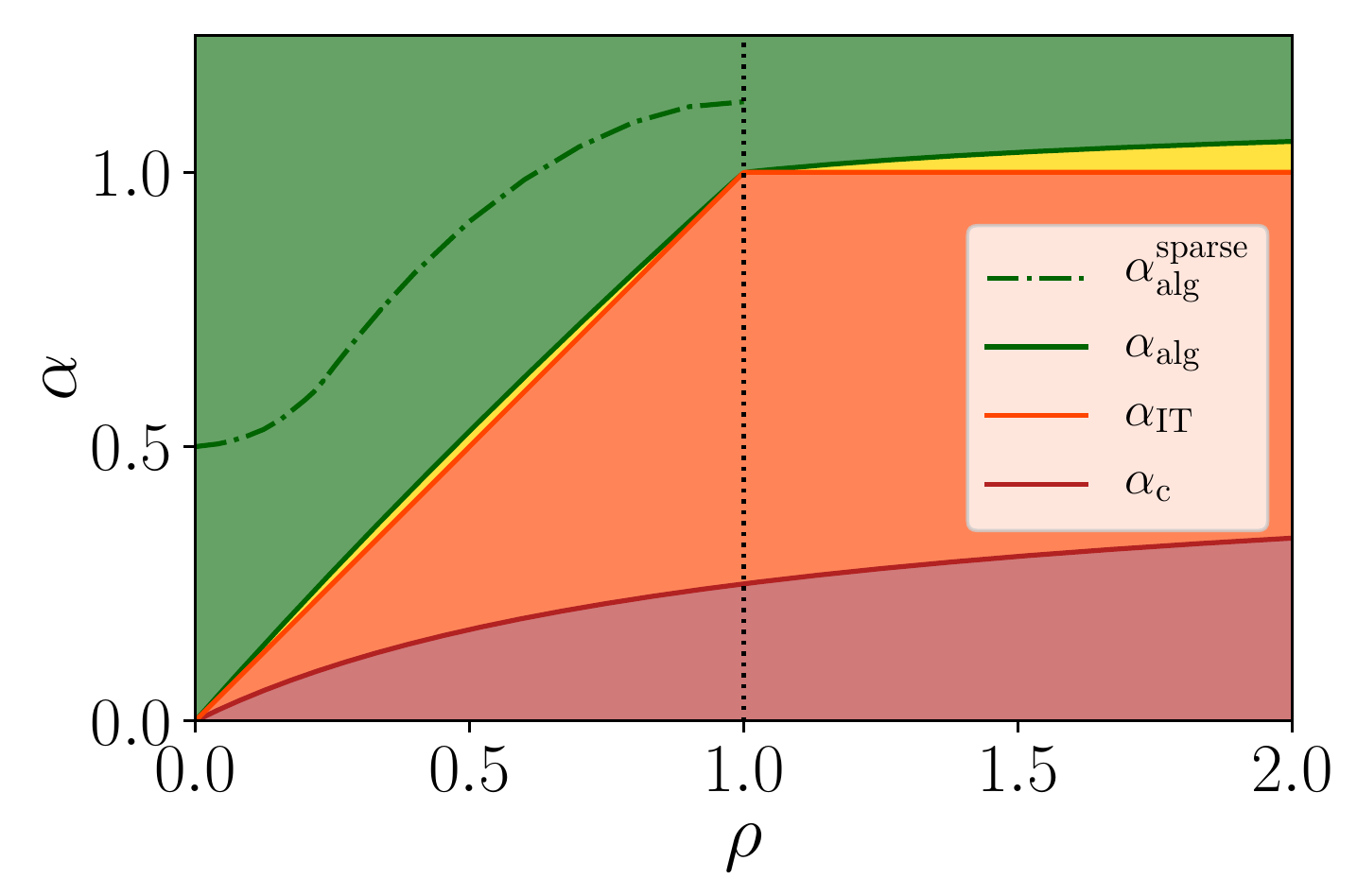}
		\includegraphics[scale=0.45]{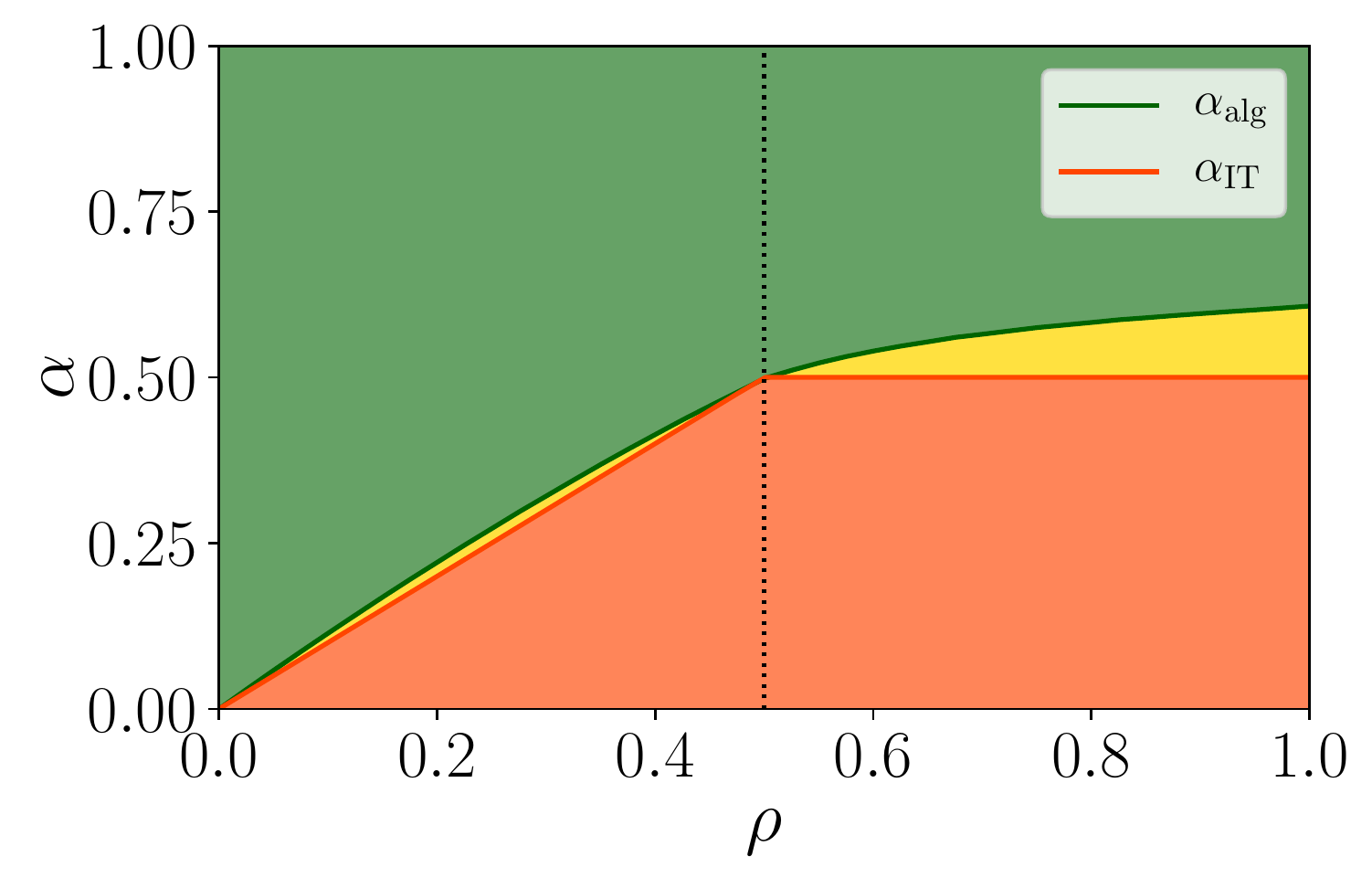}
	\caption{The same as Fig.~\ref{main:phase_diagramm_CS} for the phase retrieval problem with \textbf{(\textbf{left})} linear generative prior and \textbf{(\textbf{right})} ReLU generative prior. A major result is that while with sparse separable priors (green dashed-dotted line) compressed phase retrieval is algorithmically hard for $\alpha < 1/2$, with generative priors compressed phase retrieval is tractable down to vanishing $\alpha$ (green line). 
	In the left part we depict additionally the \textit{weak recovery} transition $\alpha_c = \rho/ [2(1+\rho)]$ (dark red line). 
	It splits the \textit{no-exact-recovery} phase  into the {\it undetectable} (dark red) and the {\it weak-recovery} region (orange). 
	}
	\label{main:phase_diagramm_PR}
\end{figure}

Fig.~\ref{main:phase_diagramm_PR} depicts the phase diagram for the phase retrieval problem with linear (left) and ReLU (right) generative priors. The information-theoretic transition is the same as the one for compressed sensing, while numerical inspection shows that $\alpha_{\rm alg}^{\rm PR} > \alpha_{\rm alg}^{\rm CS}$ for all $\rho \neq 0, 1/2,1$. In the left hand side we depict also the algorithmic transition corresponding to the sparse separable prior with non-zero components being Gaussian of zero mean, $\alpha_{\rm alg}^{\rm sparse}(\rho_s)$, as taken from \cite{Barbier2017c}. 
Crucially, in that case the algorithmic transition to exact recovery does not fall bellow $\alpha= 1/2$ even for very small (yet finite) $\rho_s$, thus effectively disabling the possibility to sense compressively. In contrast, with both the linear and ReLU generative priors we observe $\alpha_{\rm alg}(\rho \to 0) \to 0$. More specifically, the theory for the linear prior implies that $\alpha_{\rm alg}/\rho (\rho \to 0) \to \alpha_{\rm alg}^{\rm sparse}(\rho_s=1) \approx 1.128$ with the hard phase being largely reduced. 
Again the hard phase disappears entirely for $\rho=1$ for the linear model and $\rho=1/2$ for ReLU. 

\subsection{Multi-layer generative prior}
From the discussion above, we conclude that generative priors are algorithmically advantageous over sparse priors, allowing compressive sensing for the phase retrieval problem. We now investigate how the role of depth of the prior in this discussion. As before, we analyse both the linear estimation and phase retrieval problems, fixing $\sigma^{(l)} \equiv \sigma \in \{\text{linear}, \text{ReLU}\}$ at every layer $1\leq l \leq L$. Different from the $L=1$ case discussed above, for $L>1$ we have other $L-1$ free parameters characterising the layer-wise compression factors $\left(\beta_1, \dots, \beta_{L-1}\right)$.

\begin{figure}[htb!]
	\centering
		\includegraphics[scale=0.45]{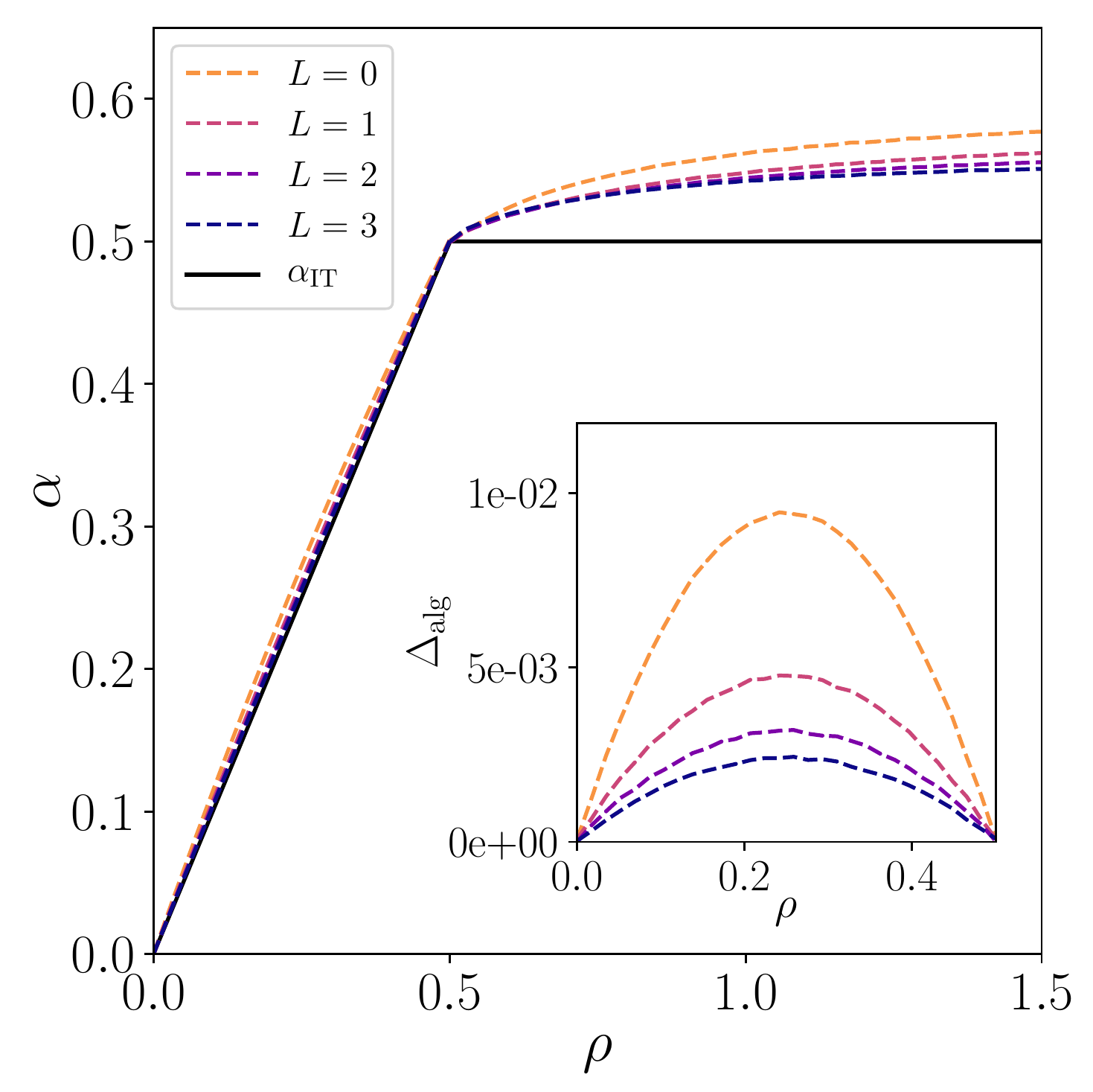}
		\includegraphics[scale=0.45]{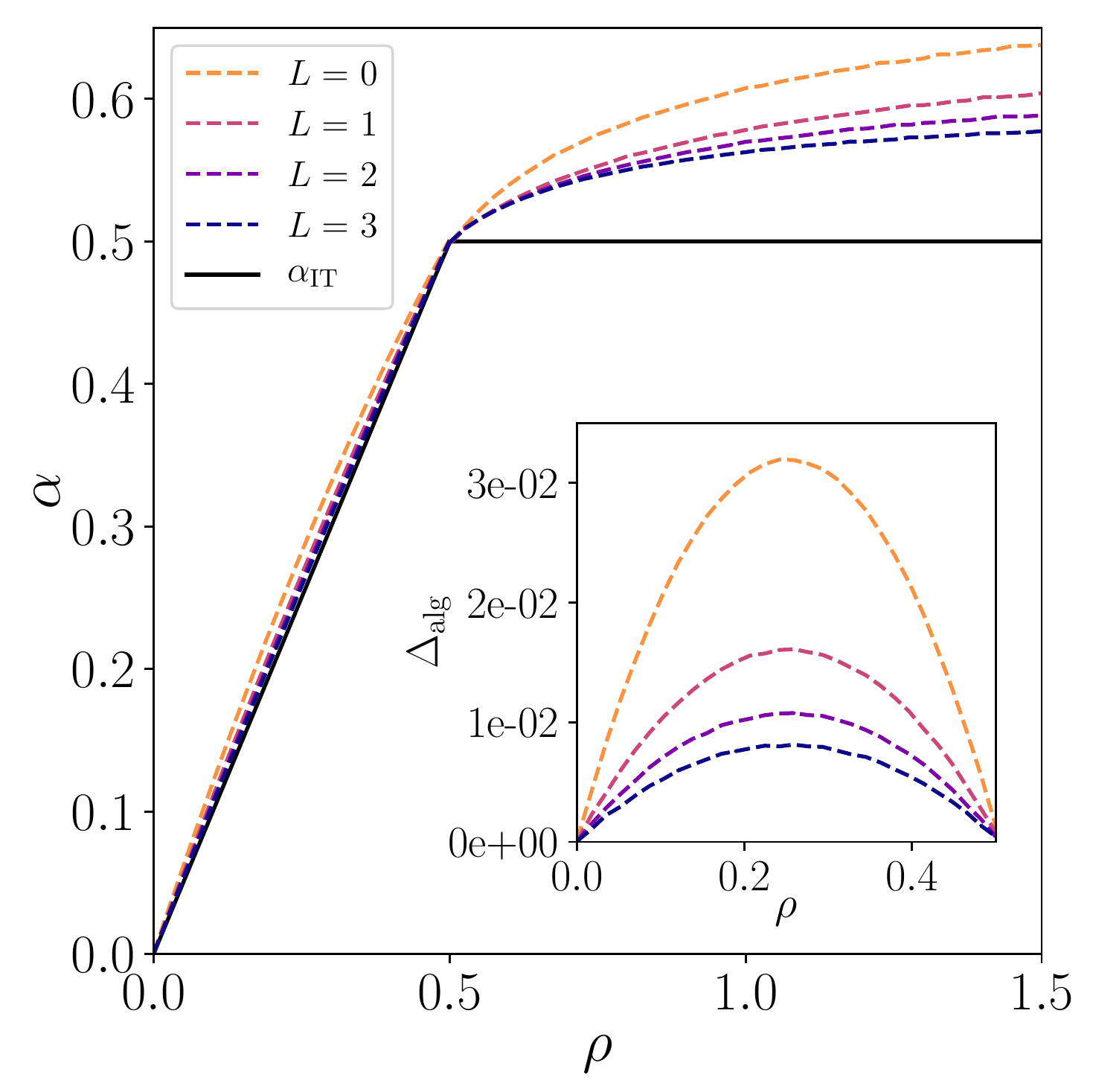}
	\caption{Phase diagrams for the compressed sensing (\textbf{left}) and phase retrieval (\textbf{right}) problems for different depths of the prior, with ReLU activation and fixed layer-wise compression $\beta_{l}=3$. Dashed lines represent the algorithmic threshold $\alpha_{\alg}$ and solid lines the perfect recovery threshold $\alpha_{\IT}$. We note that the algorithmic gap $\Delta_{\rm alg}$ (shown in insets) decreases with the network depth $L$.}
	\label{fig:multilayer:depth}
\end{figure}

First, we fix $\beta_{l}$ and investigate the role played by depth. Fig.~\ref{fig:multilayer:depth} depicts the phase diagrams for compressed sensing (left)
and phase retrieval (right) with ReLU activation with varying depth, and a fixed architecture $\beta_{l}=3$ for $1\leq l \leq L$ and note that all these curves share the same $\alpha_{\IT} = \text{min}(0.5, \rho)$. It is clear that depth improves even more the small gap already observed for a single-layer generative prior. The algorithmic advantage of multi-layer generative priors in the phase retrieval problem has been previously observed in a similar setting in \cite{hand2018phase}.

\begin{figure}[htb!]
	\centering
		\includegraphics[scale=0.45]{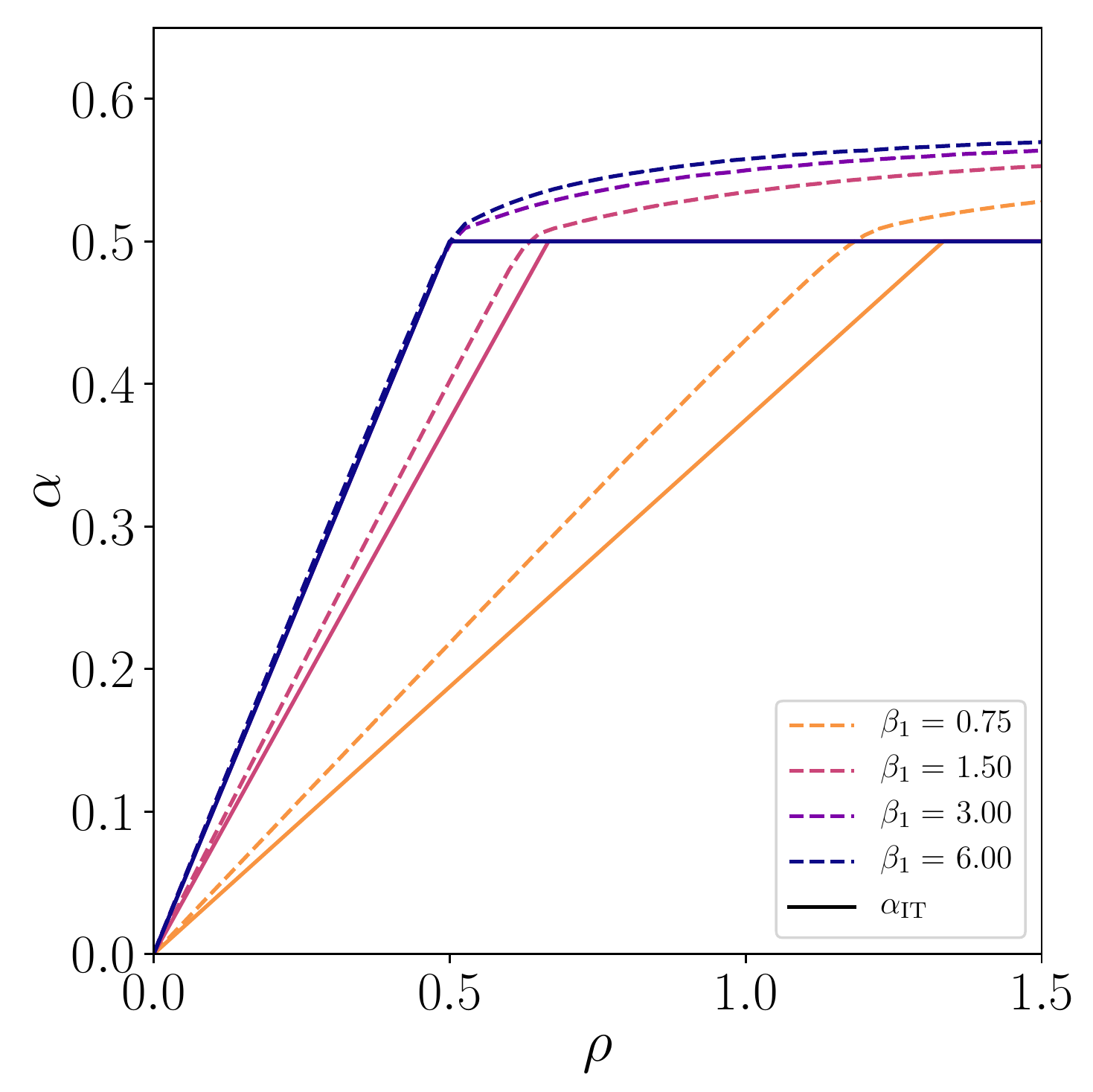}%
		\includegraphics[scale=0.45]{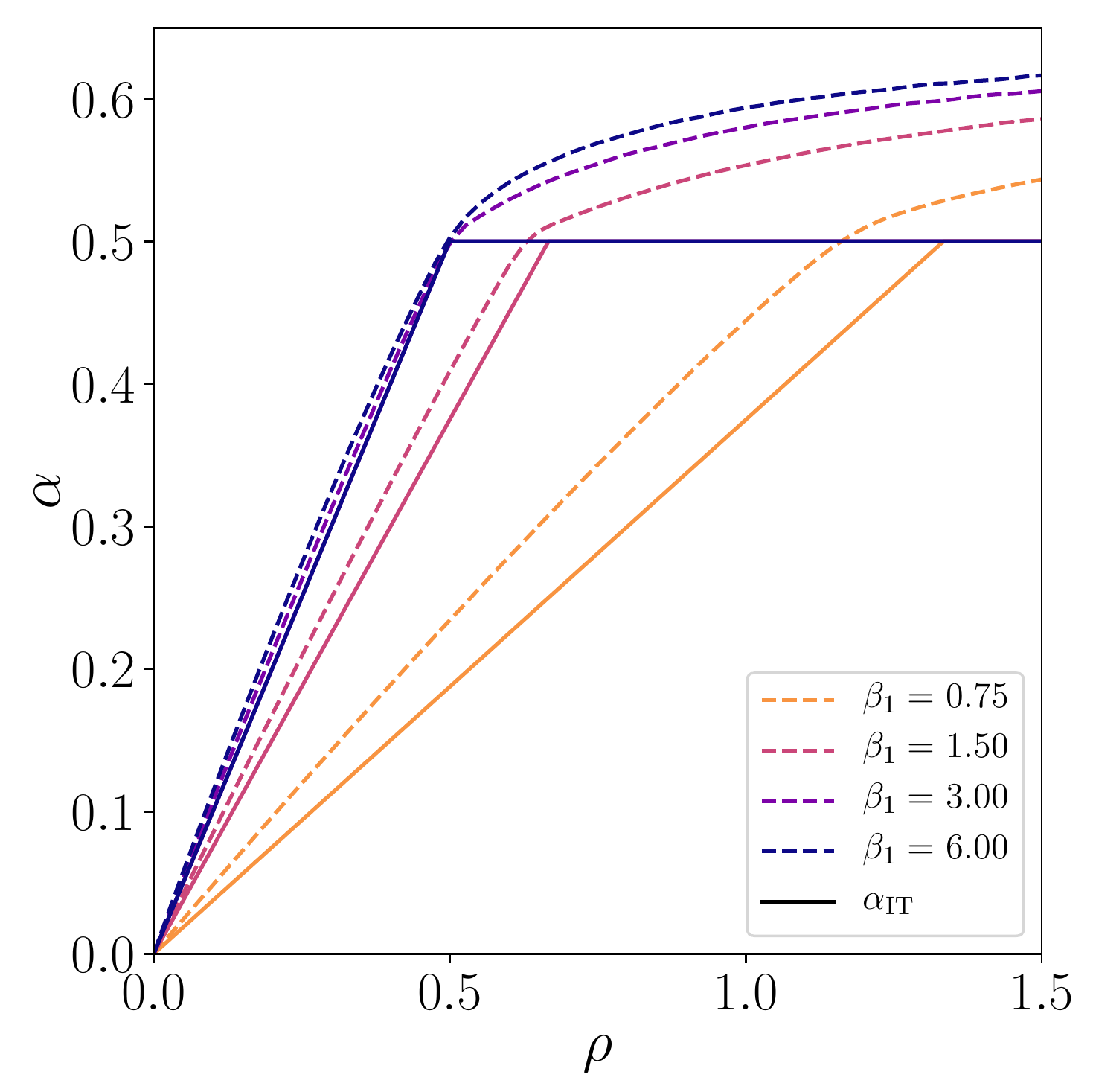}
	\caption{Phase diagrams for the compressed sensing (\textbf{left}) and phase retrieval (\textbf{right}) problems with $L=2$ and ReLU activation for different values of the layer-wise compression factor $\beta_{1}$. Dashed lines represent the algorithmic threshold $\alpha_{\alg}$ and solid lines the perfect recovery threshold $\alpha_{\IT}$. We note that for a given $\rho < 0.5$, $\alpha_{\IT}$ is decreasing with $\beta\ll 1$. However, the algorithmic gap $\Delta_{\alg}$ (shown in the inset) grows for decreasing $\beta$. Note that for $\beta_1\geq 2$ the hard phase is hardly visible at $\rho=0.5$, even though it disappears only in the large width limit, for both compressed sensing and phase retrieval settings.}	\label{fig:multilayer:compression}
	\vspace{-0.2cm}
\end{figure}

Next, we investigate the role played by the layer-wise compression factor $\beta_{l}$. Fig.~\ref{fig:multilayer:compression} depicts the phase diagrams for the compressed sensing (left) and phase retrieval (right) with ReLU activation for fixed depth $L=2$, and varying $\beta\equiv \beta_{1}$. According to the result in eq.~\eqref{eq:alphaIT}, we have $\alpha_{\IT} = \text{min}\left(1/2,\rho, 1/2\beta\right)$. It is interesting to note that there is a trade-off between compression $\beta<2$ and the algorithmic gap, in the following sense. For $\rho<0.5$ fixed, $\alpha_{\IT}$ decreases with decreasing $\beta\ll 1$: compression helps perfect recovery. However, the algorithmic gap $\Delta_{\alg}$ becomes wider for fixed $\rho<0.5$ and decreasing $\beta\ll 1$.

\begin{figure}[htb!]
	\centering
		\includegraphics[scale=0.45]{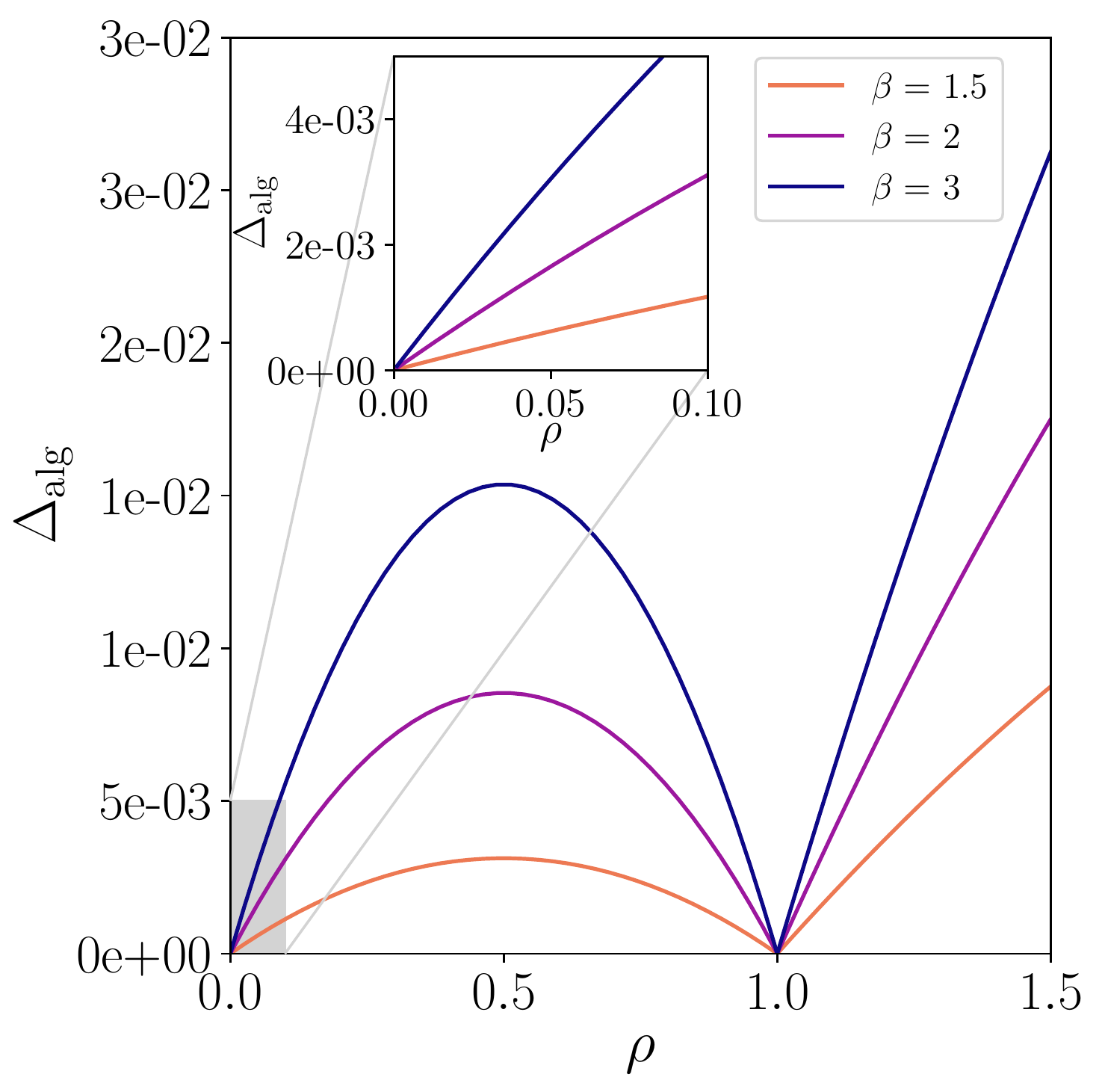}
		\includegraphics[scale=0.45]{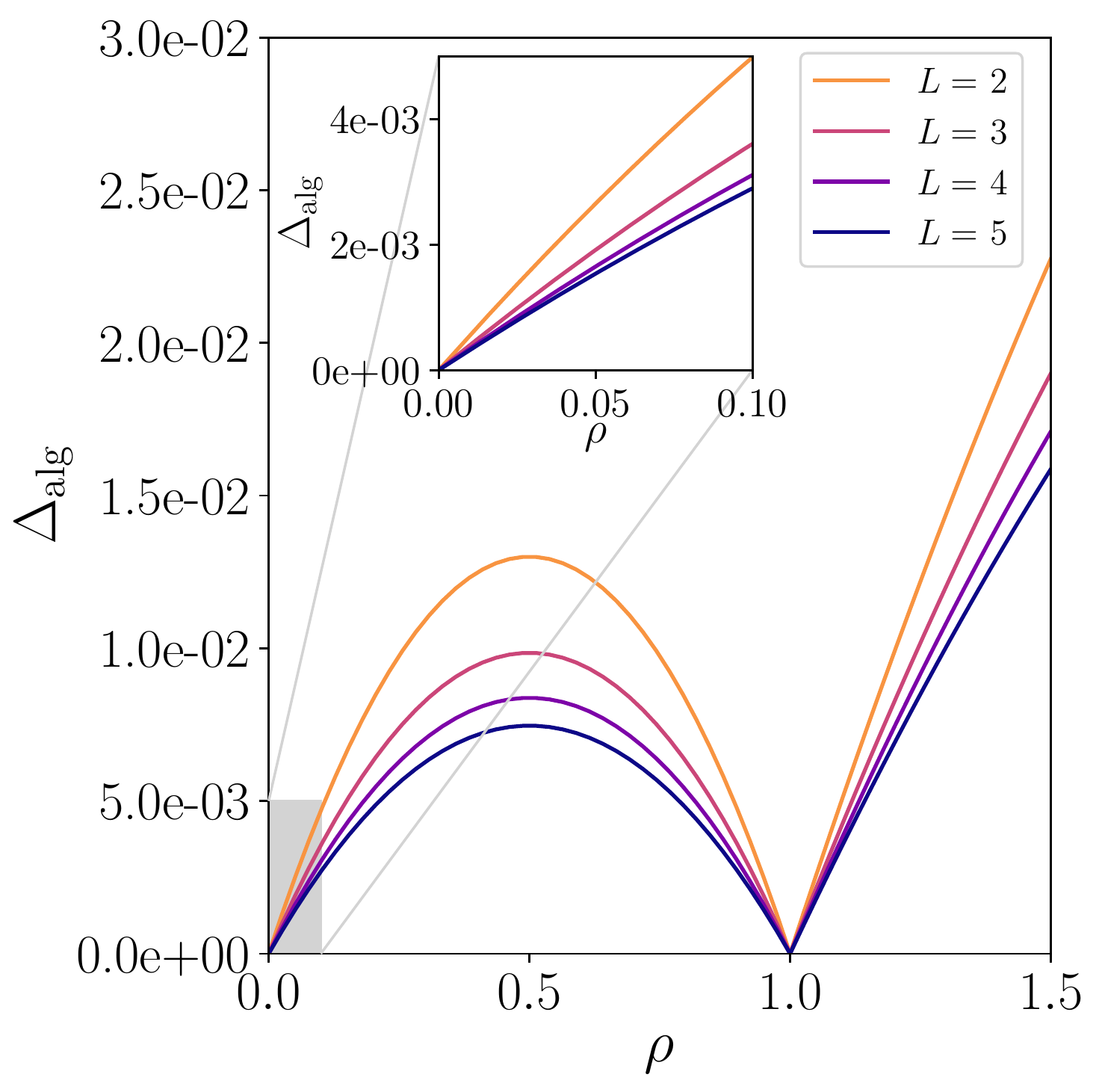}
	\caption{Algorithmic gap $\Delta_{\rm alg} = \alpha_{\rm alg} - \alpha_{\IT}$ for small $\rho$ and linear activation, as a function of (\textbf{left}) the compression $\beta\equiv \beta_{l}$ for fixed depth $L=4$ and of (\textbf{right}) depth for a fixed compression $\beta=2$.}
	\label{fig:multilayer:linear}
\end{figure}

These observations also hold for a linear generative model. In Fig.~\ref{fig:multilayer:linear} we have a closer look by plotting the algorithmic gap $\Delta_{\alg} \equiv \alpha_{\alg}-\alpha_{\IT}$ in the phase retrieval problem. On the left, we fix $L=4$ and plot the gap for increasing values of $\beta \equiv \beta_{l}$, leading to increasing $\Delta_{\alg}$. On the right, we fix $\beta = 2$ and vary the depth, observing a monotonically decreasing $\Delta_{\alg}$.

%% file: files/conclusion.tex
\section{Conclusion and perspectives}
In this manuscript we analysed how generative priors from an ensemble of random multi-layer neural networks impact signal reconstruction in the high-dimensional limit of two important inverse problems: real-valued phase retrieval and linear estimation. More specifically, we characterised the phase diagrams describing the interplay between number of measurements needed at a given signal compression $\rho$, for a range of shallow and multi-layer architectures for the generative prior. We observed that although present, the algorithmic gap significantly decreases with depth in the  studied architectures. This is particularly striking when compared with sparse priors at $\rho\ll 1$, for which the algorithmic gap is considerably wider. In practice, this means generative models given by random multi-layer neural networks allow for efficient compressive sensing in these problems.

In this work we have only considered independent random weight matrices for both the estimation layer and for the generative model. Ideally, one would like to introduce correlations in a setting closer to reality to show that the smaller computation-to-statistical gap also appears in  real-life tasks. 
The hurdle is that in those cases one does not know what is the theoretically optimal performance nor what are the optimal polynomial algorithms, so that one cannot evaluate the computation-to-statistical empirically in those cases. 
Yet another tractable case is the study of random rotationally invariant or unitary sensing matrices, as in \cite{kabashima2008inference,fletcher2018inference,barbier2018mutual,dudeja2019information}. In a different direction, it would be interesting to observe the phenomenology from this work in an experimental setting, for instance using a generative model, such as GANs or VAEs, trained on a real dataset to improve the performance of approximate message passing algorithms in a practical task.
%

%% file: Phase_retrieval_arxiv.bbl
\begin{thebibliography}{10}

\bibitem{donoho2006compressed}
David~L Donoho.
\newblock Compressed sensing.
\newblock {\em IEEE Transactions on information theory}, 52(4):1289--1306,
  2006.

\bibitem{candes2006near}
Emmanuel~J Candes and Terence Tao.
\newblock Near-optimal signal recovery from random projections: Universal
  encoding strategies?
\newblock {\em IEEE transactions on information theory}, 52(12):5406--5425,
  2006.

\bibitem{goodfellow2014generative}
Ian Goodfellow, Jean Pouget-Abadie, Mehdi Mirza, Bing Xu, David Warde-Farley,
  Sherjil Ozair, Aaron Courville, and Yoshua Bengio.
\newblock Generative adversarial nets.
\newblock In {\em Advances in neural information processing systems}, pages
  2672--2680, 2014.

\bibitem{tramel2016approximate}
Eric~W Tramel, Ang{\'e}lique Dr{\'e}meau, and Florent Krzakala.
\newblock Approximate message passing with restricted {B}oltzmann machine
  priors.
\newblock {\em Journal of Statistical Mechanics: Theory and Experiment},
  2016(7):073401, 2016.

\bibitem{tramel2016inferring}
Eric~W Tramel, Andre Manoel, Francesco Caltagirone, Marylou Gabri{\'e}, and
  Florent Krzakala.
\newblock Inferring sparsity: Compressed sensing using generalized restricted
  {B}oltzmann machines.
\newblock In {\em 2016 IEEE Information Theory Workshop (ITW)}, pages 265--269.
  IEEE, 2016.

\bibitem{bora2017compressed}
Ashish Bora, Ajil Jalal, Eric Price, and Alexandros~G Dimakis.
\newblock Compressed sensing using generative models.
\newblock In {\em Proceedings of the 34th International Conference on Machine
  Learning-Volume 70}, pages 537--546. JMLR. org, 2017.

\bibitem{manoel2017multi}
Andre Manoel, Florent Krzakala, Marc M{\'e}zard, and Lenka Zdeborov{\'a}.
\newblock Multi-layer generalized linear estimation.
\newblock In {\em 2017 IEEE International Symposium on Information Theory
  (ISIT)}, pages 2098--2102. IEEE, 2017.

\bibitem{hand2017global}
Paul Hand and Vladislav Voroninski.
\newblock Global guarantees for enforcing deep generative priors by empirical
  risk.
\newblock In {\em Conference On Learning Theory}, pages 970--978, 2018.

\bibitem{fletcher2018inference}
Alyson~K Fletcher, Sundeep Rangan, and Philip Schniter.
\newblock Inference in deep networks in high dimensions.
\newblock In {\em 2018 IEEE International Symposium on Information Theory
  (ISIT)}, pages 1884--1888. IEEE, 2018.

\bibitem{hand2018phase}
Paul Hand, Oscar Leong, and Vlad Voroninski.
\newblock Phase retrieval under a generative prior.
\newblock In {\em Advances in Neural Information Processing Systems}, pages
  9136--9146, 2018.

\bibitem{mixon2018sunlayer}
Dustin~G Mixon and Soledad Villar.
\newblock Sunlayer: Stable denoising with generative networks.
\newblock {\em arXiv preprint arXiv:1803.09319}, 2018.

\bibitem{aubin2019spiked}
Benjamin Aubin, Bruno Loureiro, Antoine Maillard, Florent Krzakala, and Lenka
  Zdeborov\'{a}.
\newblock The spiked matrix model with generative priors.
\newblock In {\em Advances in Neural Information Processing Systems 32}, pages
  8364--8375. Curran Associates, Inc., 2019.

\bibitem{reeves2017additivity}
Galen Reeves.
\newblock Additivity of information in multilayer networks via additive
  gaussian noise transforms.
\newblock In {\em 2017 55th Annual Allerton Conference on Communication,
  Control, and Computing (Allerton)}, pages 1064--1070. IEEE, 2017.

\bibitem{gabrie2018entropy}
Marylou Gabri{\'e}, Andre Manoel, Cl{\'e}ment Luneau, Nicolas Macris, Florent
  Krzakala, Lenka Zdeborov{\'a}, et~al.
\newblock Entropy and mutual information in models of deep neural networks.
\newblock In {\em Advances in Neural Information Processing Systems}, pages
  1821--1831, 2018.

\bibitem{candes2015phase}
Emmanuel~J Candes, Yonina~C Eldar, Thomas Strohmer, and Vladislav Voroninski.
\newblock Phase retrieval via matrix completion.
\newblock {\em SIAM review}, 57(2):225--251, 2015.

\bibitem{netrapalli2013phase}
Praneeth Netrapalli, Prateek Jain, and Sujay Sanghavi.
\newblock Phase retrieval using alternating minimization.
\newblock In {\em Advances in Neural Information Processing Systems}, pages
  2796--2804, 2013.

\bibitem{wu2012optimal}
Yihong Wu and Sergio Verd{\'u}.
\newblock Optimal phase transitions in compressed sensing.
\newblock {\em IEEE Transactions on Information Theory}, 58(10):6241--6263,
  2012.

\bibitem{krzakala_statistical-physics-based_2012}
F.~Krzakala, M.~M{\'e}zard, F.~Sausset, Y.~F. Sun, and L.~Zdeborov{\'a}.
\newblock Statistical physics based reconstruction in compressed sensing.
\newblock {\em Physical Review X}, 2(2), May 2012.

\bibitem{reeves2012compressed}
Galen Reeves and Michael Gastpar.
\newblock Compressed sensing phase transitions: Rigorous bounds versus replica
  predictions.
\newblock In {\em 2012 46th Annual Conference on Information Sciences and
  Systems (CISS)}, pages 1--6. IEEE, 2012.

\bibitem{zdeborova2016statistical}
Lenka Zdeborov{\'a} and Florent Krzakala.
\newblock Statistical physics of inference: Thresholds and algorithms.
\newblock {\em Advances in Physics}, 65(5):453--552, 2016.

\bibitem{barbier2016mutual}
Jean Barbier, Mohamad Dia, Nicolas Macris, and Florent Krzakala.
\newblock The mutual information in random linear estimation.
\newblock In {\em 2016 54th Annual Allerton Conference on Communication,
  Control, and Computing (Allerton)}, pages 625--632. IEEE, 2016.

\bibitem{Reeves}
G.~Reeves and H.~D. Pfister.
\newblock The replica-symmetric prediction for compressed sensing with gaussian
  matrices is exact.
\newblock In {\em 2016 IEEE International Symposium on Information Theory
  (ISIT)}, pages 665--669, July 2016.

\bibitem{Barbier2017c}
Jean Barbier, Florent Krzakala, Nicolas Macris, L{\'e}o Miolane, and Lenka
  Zdeborov{\'a}.
\newblock Optimal errors and phase transitions in high-dimensional generalized
  linear models.
\newblock {\em Proceedings of the National Academy of Sciences},
  116(12):5451--5460, 2019.

\bibitem{mezard1987spin}
Marc M{\'e}zard, Giorgio Parisi, and Miguel Virasoro.
\newblock {\em Spin glass theory and beyond: An Introduction to the Replica
  Method and Its Applications}, volume~9.
\newblock World Scientific Publishing Company, 1987.

\bibitem{donoho2009message}
David~L Donoho, Arian Maleki, and Andrea Montanari.
\newblock Message-passing algorithms for compressed sensing.
\newblock {\em Proceedings of the National Academy of Sciences},
  106(45):18914--18919, 2009.

\bibitem{rangan2011generalized}
Sundeep Rangan.
\newblock Generalized approximate message passing for estimation with random
  linear mixing.
\newblock In {\em 2011 IEEE International Symposium on Information Theory
  Proceedings}, pages 2168--2172. IEEE, 2011.

\bibitem{schniter2014compressive}
Philip Schniter and Sundeep Rangan.
\newblock Compressive phase retrieval via generalized approximate message
  passing.
\newblock {\em IEEE Transactions on Signal Processing}, 63(4):1043--1055, 2014.

\bibitem{metzler2017coherent}
Christopher~A Metzler, Manoj~K Sharma, Sudarshan Nagesh, Richard~G Baraniuk,
  Oliver Cossairt, and Ashok Veeraraghavan.
\newblock Coherent inverse scattering via transmission matrices: Efficient
  phase retrieval algorithms and a public dataset.
\newblock In {\em 2017 IEEE International Conference on Computational
  Photography (ICCP)}, pages 1--16. IEEE, 2017.

\bibitem{bayati2011dynamics}
Mohsen Bayati and Andrea Montanari.
\newblock The dynamics of message passing on dense graphs, with applications to
  compressed sensing.
\newblock {\em IEEE Transactions on Information Theory}, 57(2):764--785, 2011.

\bibitem{krzakala2012probabilistic}
Florent Krzakala, Marc M{\'e}zard, Francois Sausset, Yifan Sun, and Lenka
  Zdeborov{\'a}.
\newblock Probabilistic reconstruction in compressed sensing: algorithms, phase
  diagrams, and threshold achieving matrices.
\newblock {\em Journal of Statistical Mechanics: Theory and Experiment},
  2012(08):P08009, 2012.

\bibitem{tao2009}
Terence Tao.
\newblock Compressed sensing, or: the equation ${A}x=b$, revisited.
\newblock Clay-Mahler Lecture Series, 2009.

\bibitem{baker2020tramp}
Antoine Baker, Benjamin Aubin, Florent Krzakala, and Lenka Zdeborov{\'a}.
\newblock Tramp: Compositional inference with tree approximate message passing.
\newblock {\em arXiv preprint arXiv:2004.01571}, 2020.

\bibitem{kabashima2008inference}
Yoshiyuki Kabashima.
\newblock Inference from correlated patterns: a unified theory for perceptron
  learning and linear vector channels.
\newblock In {\em Journal of Physics: Conference Series}, volume~95, page
  012001. IOP Publishing, 2008.

\bibitem{barbier2018mutual}
Jean Barbier, Nicolas Macris, Antoine Maillard, and Florent Krzakala.
\newblock The mutual information in random linear estimation beyond iid
  matrices.
\newblock In {\em 2018 IEEE International Symposium on Information Theory
  (ISIT)}, pages 1390--1394. IEEE, 2018.

\bibitem{dudeja2019information}
Rishabh Dudeja, Junjie Ma, and Arian Maleki.
\newblock Information theoretic limits for phase retrieval with subsampled
  {H}aar sensing matrices.
\newblock {\em arXiv preprint arXiv:1910.11849}, 2019.

\end{thebibliography}
